\documentclass[a4paper, 11pt, oneside]{amsart}
\usepackage{amsmath, amsthm,amssymb, amscd}
\usepackage{ifthen}

 \provideboolean{showlabels}
 \setboolean{showlabels}{false}
\newcommand{\mylabel}[1]{\ifthenelse{\boolean{showlabels}}{{\tt{[{#1}]}}
\label{#1}}{\label{#1}}}

\theoremstyle{plain}
\newtheorem{theorem}[equation]{Theorem}
\newtheorem{lemma}[equation]{Lemma}
\newtheorem{proposition}[equation]{Proposition}
\newtheorem{corollary}[equation]{Corollary}
\newtheorem*{main}{Main Theorem}
\newtheorem*{A}{Theorem A}
\newtheorem{step}{Step}

\theoremstyle{definition}

\numberwithin{equation}{section}

\DeclareMathOperator{\Core}{Core}
\DeclareMathOperator{\Dom}{Dom}

\DeclareMathOperator{\Irr}{Irr}
\DeclareMathOperator{\Ker}{Ker}
\DeclareMathOperator{\Lin}{Lin}
\DeclareMathOperator{\rad}{rad}

\DeclareMathOperator{\C}{C}

\DeclareMathOperator{\N}{N}
\DeclareMathOperator{\SSS}{\bar{\N}}
\DeclareMathOperator{\Z}{Z}
\DeclareMathOperator{\GL}{GL}
\DeclareMathOperator{\LL}{LL}
\DeclareMathOperator{\LR}{LR}

\DeclareMathOperator{\ST}{ST}
\DeclareMathOperator{\ZG}{ZG}

\DeclareMathOperator{\MLR}{MLR}
\DeclareMathOperator{\ZST}{ZST}

\newcommand{\CC}{\mathbb{C}}

\newcommand{\LSG}{\mathcal{L}}

\newcommand{\tin}{\text{\textrm{\textup{ in }}}}

\newcommand{\mfrt}{\mathfrak{T}}

\newcommand{\mfrl}{\mathfrak{L}}

\newcommand{\zetaQ}[1]{\zeta^{(#1)}}

\newcommand{\normaleq}{\trianglelefteq}
\newcommand{\normal}{\triangleleft}
\newcommand{\sq}[1]{\tilde{#1}}

\newcommand{\gen}[1]{\left < #1 \right >}

\begin{document}

\title[Dade and Loukaki, Linear Limits]{ Linear limits of irreducible characters }

\author{Everett Dade} 
\address{Department of Mathematics\\
University of Illinois at Urbana-Champaign\\
1409 W.~Green St.\\
Urbana, IL 61801, USA}
\email{dade@math.uiuc.edu}

\author{Maria Loukaki}
\address{School of Mathematics\\
Georgia Institute of Technology\\
686 Cherry Street NW\\
Atlanta, GA 30332-0160, USA}
\email{loukaki@math.gatech.edu}

%\abstract{This is the abstract}

\date{December 13, 2004} 

\thanks{ The second author was 
 partially supported by NSF grants DMS 96-00106 and DMS 99-70030} 

\maketitle

\setcounter{section}{0}

\textbf{Abstract.}  Nearly twenty years ago Isaacs and the first author
of this paper wrote a series of articles \cite{isa2}, \cite{da3},
\cite{da2} about what were called ``stabilizer limits'' of group
characters, following the terminology of Berger \cite{be}.  The second
author, in her thesis \cite{lo}, needed one of the results of those
articles in a new situation which was not treated earlier.  Eventually
she was able, by complicated and delicate arguments, to reduce her proof
to a special case where \cite[Theorem 8.4]{da3} could be applied.  But
this approach was extremely awkward.  In the present paper we use
arguments similar to those in the earlier articles to prove a Main
Theorem from which the exact Theorem A needed in \cite{lo} easily
follows.

\textbf{Introduction.} To explain the precise content of our theorems
requires a lot of notation and terminology.  We start with the family
$\mfrt$ of all ordered triples $(G,N,\psi)$, consisting of a finite
group $G$, a normal subgroup $N$ of $G$, and an irreducible complex
character $\psi$ of $N$. Fix a triple $T = (G,N,\psi)$ in $\mfrt$. This
$T$ determines several other objects, such as the stabilizer $G(\psi)$
of the character $\psi$ in the group $G$, the image $G(\psi)/N$ of that
stabilizer in the factor group $G/N$, and the character $\psi^G$ of $G$
induced by $\psi$. We define the \emph{center} $\Z(T)$ of $T$ to be the
center $\Z(\psi^G)$ of this induced character, and the \emph{central
character} of $T$ to be the unique linear character $\zetaQ{T}$ of
$\Z(T)$ such that $\psi^G$ restricts to a multiple of $\zetaQ{T}$ on
$\Z(T)$. Then $\Z(T)$ is a normal subgroup of $G$ contained in $N$, and
$\zetaQ{T}$ is a $G$-invariant linear character of $\Z(T)$ lying under
$\psi$.

By a \emph{subtriple} of $T$ we mean any triple $T' = (G',N',\psi')$,
consisting of a subgroup $G'$ of $G$, the intersection $N' = G' \cap N$ of
that subgroup with $N$, and an irreducible complex character $\psi'$ of
$N'$ lying under the irreducible character $\psi$ of $N$. Any such $T'$
also lies in $\mfrt$.  Furthermore, inclusion $G' \hookrightarrow G$
induces a natural embedding of the factor group $G'/N'$ into $G/N$. It is
easy to see that the intersection $G' \cap \Z(T)$ is contained in
$\Z(T')$, and that both $\zetaQ{T}$ and $\zetaQ{T'}$ restrict to the same
linear character of this intersection. In general there is no useful
relation between the stabilizers $G'(\psi')$ and $G(\psi)$, or even between
their images $G'(\psi')/N'$ and $G(\psi)/N$. Notice that any subtriple
$T''$ of the subtriple $T'$ is also a subtriple of $T$. So is any
\emph{conjugate} $(T')^{\tau} = ((G')^{\tau}, (N')^{\tau}, (\psi')^{\tau})$
of $T'$ by any $\tau \in G(\psi)$.

Fix a subtriple $T' = (G',N',\psi')$ of $T$. We say that $T'$ is a
\emph{linear reduction} of $T$ if there exist some normal subgroup $L$ of
$G$ contained in $N$, and some complex linear character
$\lambda$ of $L$ lying under $\psi$, such that $G'$ is the stabilizer
$G(\lambda)$ of $\lambda$ in $G$, and $\psi'$ is the unique irreducible
character of the stabilizer $N' = G' \cap N = N(\lambda)$ lying over
$\lambda$ and inducing $\psi$.  We say that $T'$ is a \emph{multilinear
reduction} of $T$ if there is some finite chain $T_0, T_1, \dots, T_n$ of
subtriples of $T$, starting with $T_0 = T$ and ending with $T_n = T'$, such
that each $T_i$, for $i = 1,2,\dots,n$, is a linear reduction of its
predecessor $T_{i-1}$. In that case we call $T_0,T_1,\dots,T_n$ a \emph{linear reducing chain} from $T$ to $T'$. We say that $T'$ is a \emph{linear limit} of $T$ if
it is a multilinear reduction of $T$ such that the only possible linear
reduction of $T'$ is $T'$ itself. Each of the classes of linear reductions,
of multilinear reductions, or of linear limits of $T$ is closed under
conjugation by elements of $N$. If $T'$ belongs to any of these classes,
then the natural embedding of $G'/N'$ into $G/N$  sends $G'(\psi')/N'$
isomorphically onto $G(\psi)/N$, the center $\Z(T)$ is a subgroup of
$\Z(T')$, and the central character $\zetaQ{T}$ is the restriction of
$\zetaQ{T'}$ to $\Z(T)$.

We say that the subtriple $T'$ \emph{covers $T$ modulo $\Z(T)$} if the
subgroup $G'\Z(T)$ is equal to $G$. In this case the natural embedding is
an isomorphism of the factor group $G'/N'$ onto $G/N$. Furthermore, the
intersection $G' \cap \Z(T)$ is exactly $\Z(T')$, so that inclusion $G'
\hookrightarrow G$ induces an isomorphism of $G'/\Z(T')$ onto $G/\Z(T)$.
The character $\psi'$ is now the restriction of $\psi$ to $N'$, and the
stabilizer $G(\psi)$ is the product $G'(\psi')\Z(T)$.  It follows that the
natural isomorphism of $G'/N'$ onto $G/N$ sends the subgroup $G'(\psi')/N'$
onto $G(\psi)/N$.

After all these definitions we can finally state our
\begin{main} If\/ $T'$ and $T''$ are two linear limits of some triple $T
\in \mfrt$, then there exists a finite chain $T_0, T_1, \dots, T_n$ of
subtriples $T_i$ of\/ $T$, starting with $T_0 = T'$ and ending with $T_n =
T''$, such that, for each $i = 1,2,\dots,n$, either $T_i$ is an
$N$-conjugate of\/ $T_{i-1}$, or $T_i$ covers $T_{i-1}$ modulo
$\Z(T_{i-1})$, or $T_{i-1}$ covers $T_i$ module $\Z(T_i)$.
\end{main}

Suppose that the factor group $\SSS(T') = N'/\Z(T')$ is nilpotent, for
some linear limit $T' = (G',N',\psi')$ of $T$. Then the Hall-Higman
argument tells us that $\SSS(T')$ is abelian, and that commutation $\sigma,
\tau \mapsto [\sigma,\tau] = \sigma^{-1}\tau^{-1}\sigma\tau$ in $N'$
induces a non-singular, strongly alternating, $G'$-invariant, bilinear form
$\bar c'$ from $\SSS(T') \times \SSS(T')$ to the multiplicative group
$\CC^{\times}$ of the complex numbers, sending $\sigma\Z(T'), \tau\Z(T')$
to $\zetaQ{T'}([\sigma,\tau])$, for any $\sigma, \tau \in N'$. This
bilinear form, and the natural  action of the factor group $G'/N'$ on the
abelian group $\SSS(T') = N'/\Z(T')$, turn $\SSS(T')$ into a
\emph{symplectic $G'/N'$-group}. 

Since $\SSS(T')$ is symplectic, the irreducible character $\psi'$ of $N'$
lying over $\zetaQ{T'}$ is zero outside $\Z(T')$, and a multiple of
$\zetaQ{T'}$ on $\Z(T')$. Because $\zetaQ{T'}$ is $G'$-invariant, so is
$\psi'$. Hence the natural embedding sends $G'/N' = G'(\psi')/N'$
isomorphically onto $G(\psi)/N$. We use this isomorphism to transfer the
action of $G'/N'$ on $\SSS(T')$ to one of $G(\psi)/N$ on that group. In
this way $\SSS(T')$ becomes a symplectic $G(\psi)/N$-group

The above Main Theorem implies that the factor group $\SSS(T'')
=N''/\Z(T'')$, for any other linear limit $T'' = (G'', N'', \psi'')$ of
$T$, is isomorphic to $\SSS(T')$, and hence is nilpotent. So $\SSS(T')$ is
also a symplectic $G(\psi)/N$-group. For us the most important
consequence of the Main Theorem is
\begin{A} If some linear limit $T' = (G',N',\psi')$ of a triple $T  =
(G,N,\psi) \in \mfrt$ has a nilpotent factor group $\SSS(T') = N'/\Z(T')$,
then $\SSS(T')$ is naturally a symplectic $G(\psi)/N$-group. So
is\/ $\SSS(T'') = N''/\Z(T'')$, for any other linear limit $T'' =
(G'',N'',\psi'')$ of\/ $T$. Furthermore, $\SSS(T')$ and $\SSS(T'')$ are
isomorphic as symplectic $G(\psi)/N$-groups.
\end{A}

\noindent
This is exactly the theorem needed in \cite{lo}.

The first section of this paper gives most of the notation we use for groups and characters, along with a few specialized remarks about centers and central products of characters.  The second, third, and fourth sections discuss, respectively, the triples in $\mfrt$, their linear reductions, and their multilinear reductions.  Some of the results about multilinear reductions at the end of the fourth section are rather startling. For example, there is actually a canonical linear reducing chain from any $T \in \mfrt$ to any given multilinear reduction $T'$ of $T$ (see the remarks following Proposition \ref{CanonLProp}).  Furthermore, Proposition \ref{MLRSubProp} says that any such $T'$ is also a multilinear reduction of any subtriple $\sq T$ of $T$ such that $T'$ is a subtriple of $\sq T$.  

The fifth section below is a detailed study of covering triples and the related correspondences of subgroups, characters, triples, linear and multilinear reductions, and linear limits. The above Main Theorem is proved as Theorem \ref{LLThm} in the sixth section, while Theorem A  is proved as Theorem \ref{AThm} in the seventh and final section of the paper.

\section{ Groups and characters  } \mylabel{GpsNChars}

Unless otherwise specified, any group we use will be multiplicative. We
write the indentity element of a group $G$ as $1 = 1_G$, the center of $G$
as $\Z(G)$, the conjugate of an element $\sigma \in G$ by an element $\tau
\in G$ as $\sigma^{\tau} = \tau^{-1}\sigma\tau$, and the commutator of
$\sigma$ with $\tau$ as $[\sigma, \tau] = \sigma^{-1}\tau^{-1}\sigma\tau$.
The expressions $H \subseteq G$, $H \le G$ or $H \normaleq G$ mean that $H$
is, respectively, a subset, a subgroup or a normal subgroup of $G$. To say,
in addition, that $H$ is properly contained in $G$, we write $H \subset G$,
$H < G$ or $H \normal G$, respectively. If $H, K \le G$, then the
normalizer of $H$ in $K$ is denoted by either $\N_K(H)$ or $\N(H \tin K)$,
and the centralizer of $H$ in $K$ by either $\C_K(H)$ or $\C(H \tin K)$. We
write $\Core_K(H)$ for the $K$-core
\begin{equation} \mylabel{Core} \Core_K(H) = \bigcap_{\tau \in K}H^{\tau}
\end{equation}
of $H$, the largest subgroup of $H$ normalized by $K$. As usual, the
``commutator'' $[H,K]$ is the subgroup of $G$ generated by
all the commutators $[\sigma, \tau]$ for $\sigma \in H$ and $\tau \in K$,
and not the set of those commutators.  We always write the commutator
subgroup of $G$ as $[G,G]$, reserving the notation $G'$ for other uses.

For the rest of this section $G$ will be a fixed, but arbitrary, group of
finite order. We write that order as $|G|$, and the index of any subgroup
$H \le G$ in $G$ as $[G : H]$. By a \emph{character} $\psi$ of $G$ we
always mean a non-zero complex character, considered as a function from $G$
to the field $\CC$ of all complex numbers. Such a  $\psi$ determines $G$ as
its domain $\Dom(\psi)$ of definition.  We write $\gen{\psi,\phi} =
\gen{\psi,\phi}_G$ for the usual inner product of two characters $\psi,
\phi$ of $G$. We denote by $\Irr(G)$ the finite
set of all irreducible characters of $G$, and by $\hat 1 = \hat 1_G$ the
trivial character, sending each $\sigma \in G$ to $1 = 1_{\CC}$.  Then
$\hat 1$ is the identity element in the finite group $\Lin(G)$ of all
linear characters of $G$, i.~e., of all homomorphisms of $G$ into the
multiplicative group $\CC^{\times}$ of the complex field $\CC$. Of
course, $\Lin(G)$ is also a subset of $\Irr(G)$.

If $H$ is a subgroup of $G$, then we write $\psi_H$ for the character of
$H$ restricted from a given character $\psi$ of $G$, and $\phi^G$ for the
character of $G$ induced by a given character $\phi$ of $H$. Conjugation
by any $\tau \in G$ sends $H$ to the conjugate subgroup $H^{\tau} =
\tau^{-1}H\tau$ of $G$. It also sends any character $\phi$ of $H$ to the
conjugate character $\phi^{\tau}$ of $H^{\tau}$, defined by
\begin{subequations} \mylabel{ConjCharEqs}
\begin{equation} \mylabel{ConjChar} \phi^{\tau}(\sigma^{\tau}) =
\phi(\sigma) \in \CC
\end{equation}
for any $\sigma \in H$. In this way the group $G$ acts by conjugation on
the family of all characters of subgroups of $G$. We denote by $K(\phi)$
the stabilizer
\begin{equation} \mylabel{CharStab}
K(\phi) = \{ \tau \in K \mid \phi^{\tau} = \phi\,\}
\end{equation}
of any character $\phi$ of $H$ in any subgroup $K \le G$. Notice that
$K(\phi)$ is contained in $\N_K(H)$, since any $\tau \in K$ fixing $\phi$
must normalize $H = \Dom(\phi)$. More generally, we write $K(\phi_1,
\phi_2, \dots,\phi_m)$ for the common stabilizer
\begin{equation} \mylabel{CharsStab} K(\phi_1, \phi_2, \dots, \phi_m) =
K(\phi_1) \cap K(\phi_2) \cap \dots\cap K(\phi_m)
\end{equation}
\end{subequations}
in $K$ of any characters $\phi_1, \phi_2, \dots, \phi_m$ of subgroups
$H_1, H_2, \dots, H_m$, respectively, of $G$.

We're going to use several properties of the kernel and center of an
arbitrary character $\psi$ of $G$. For the benefit of the reader (and
perhaps of ourselves) we recall those properties here. We first choose
some complex matrix representation $A \colon \sigma \mapsto A(\sigma)$ of
$G$ affording the character $\psi$. Then $A$ is a homomorphism of the group
$G$ into the group $\GL_n(\CC)$ of all non-singular $n \times
n$ complex matrices, where $n$ is the degree $\psi(1)$ of $\psi$.  The
\emph{kernel} $\Ker(\psi)$ of $\psi$ is just the kernel of this
homomorphism, the normal subgroup of $G$ consisting of all $\sigma \in G$
such that $A(\sigma)$ is the identity matrix $I \in \GL_n(\CC)$. The
\emph{center} $\Z(\psi)$ of $\psi$ is the larger normal subgroup of $G$
consisting of all $\sigma \in G$ such that $A(\sigma)$ is the product
$\zetaQ{\psi}(\sigma)I$ for some $\zetaQ{\psi}(\sigma) \in \CC^{\times}$.
The resulting function $\zetaQ{\psi} \colon \Z(\psi) \to \CC^{\times}$ is a
$G$-invariant linear character of $\Z(\psi)$ lying under $\psi$. We call
$\zetaQ{\psi}$ the \emph{central character} for $\psi$. Of course, the
objects $\Ker(\psi)$, $\Z(\psi)$ and $\zetaQ{\psi}$ depend only on $\psi$,
and not on the choice of the representation $A$ affording $\psi$.

One obvious consequence of the above definitions is that $\zetaQ{\psi}$ and
$\psi$ have the same kernel
\begin{equation} \mylabel{KerZetaPsi} \Ker(\zetaQ{\psi}) = \Ker(\psi).
\end{equation}
It follows that the factor group $\Z(\psi)/\Ker(\psi)$ is cyclic
as well as central in $G/\Ker(\psi)$. Another immediate consequence is that
\begin{equation} \mylabel{PsiZetaPsi} \psi(\sigma\tau) = \psi(\tau\sigma) =
\psi(\sigma)\zetaQ{\psi}(\tau) \in \CC
\end{equation}
for any $\sigma \in G$ and $\tau \in \Z(\psi)$. Since each matrix
$A(\sigma)$, for $\sigma \in G$, is diagonalizable, a third useful
consequence can be stated as
\begin{proposition} \mylabel{UnderZetaPsiProp} A subgroup $L$ of\/ $G$ is
contained in $\Z(\psi)$ if and only if the restriction $\psi_L$ is a
multiple $\psi(1)\lambda$ of a single linear character $\lambda \in
\Lin(L)$. In that case $\lambda$ is the restriction $(\zetaQ{\psi})_L$ of\/
$\zetaQ{\psi}$. Hence $\Ker{\lambda} \le \Ker(\zetaQ{\psi}) = \Ker(\psi)$. 
\end{proposition}

We now apply the above considerations when $\psi$ is the character
$\phi^G$ of $G$ induced by some character $\phi$ of some subgroup $H \le
G$. We may choose the representation $A$ affording $\psi$ to be induced
from a matrix representation $B$ of $H$ affording $\phi$.  Then $A$ sends
an element $\sigma \in G$ to $z$ times an identity matrix, for some $z \in
\CC^{\times}$, if and only if right multiplication by $\sigma$ sends each
coset $H\tau$, for $\tau \in G$, onto itself $H\tau\sigma = H\tau$, and
$B$ sends the corresponding element $\tau\sigma\tau^{-1} \in H$ to $z$
times an identity matrix. This just says that $\Z(\phi^G)$ consists of all
elements $\sigma$ in the $G$-core $\Core_G(\Z(\phi)) = \bigcap_{\tau \in
G}\Z(\phi)^{\tau}$ of $\Z(\phi)$ such that $(\zetaQ{\phi})^{\tau}(\sigma)$
has the same value $z$ for all $\tau \in G$. Then
$\zetaQ{\phi^G}(\sigma)$ is that common value $z$. From this discussion we
conclude that $\Ker(\phi^G)$ is the $G$-core
\begin{subequations} \mylabel{PhiGEqs}
\begin{equation} \mylabel{KerPhiG} \Ker(\phi^G) = \Core_G(\Ker(\phi)) =
\bigcap_{\tau \in G} \Ker(\phi)^{\tau} = \bigcap_{\tau \in G}
\Ker(\phi^{\tau}) \le \Ker(\phi)
\end{equation}
of $\Ker(\phi)$, that $\Z(\phi^G)$ is a subgroup
\begin{equation} \mylabel{ZPhiGInZPhi} \Z(\phi^G) \le \Z(\phi)
\end{equation}
of $\Z(\phi)$, and that $\zetaQ{\phi^G}$ is the restriction
\begin{equation} \mylabel{ZetaPhiG} \zetaQ{\phi^G} =
(\zetaQ{\phi})_{\Z(\phi^G)} 
\end{equation}
\end{subequations}
of $\zetaQ{\phi} \in \Lin(\Z(\phi))$ to that subgroup. Another consequence
is
\begin{proposition} \mylabel{UnderZetaPhiGProp} The center $\Z(\phi^G)$ is
the largest normal subgroup $L$ of\/ $G$ contained in $H$ such that
$\phi_L$ is a multiple $\phi(1)\lambda$ of some $G$-invariant $\lambda \in
\Lin(L)$. Any such $\lambda$ is the restriction $(\zetaQ{\phi^G})_L$ of\/
$\zetaQ{\phi^G}$ to the corresponding $L$.  Hence $\Ker(\lambda) \le
\Ker(\zetaQ{\phi^G})  = \Ker(\phi^G) \le \Ker(\phi)$.
\end{proposition}
\begin{proof} Since $\phi^G$ is a character of $G$, its center $\Z(\phi^G)$
is a normal subgroup of $G$, and its central character $\zetaQ{\phi^G} \in
\Lin(\Z(\phi^G))$ is $G$-invariant. We know from (\ref{PhiGEqs}b,c) that
$\zetaQ{\phi^G}$ is the restriction of $\zetaQ{\phi}$ to the subgroup
$\Z(\phi^G)$ of $\Z(\phi) \le H$. So
\[ \phi_{\Z(\phi^G)} = (\phi_{\Z(\phi)})_{\Z(\phi^G)} =
(\phi(1)\zetaQ{\phi})_{\Z(\phi^G)} = \phi(1)\zetaQ{\phi^G}. \]

Suppose that $L$ is any normal subgroup of $G$ contained in $H$ such that
$\phi$ restricts to a multiple $\phi_L = \phi(1)\lambda$ of some
$G$-invariant $\lambda \in \Lin(L)$.  Then $L = L^{\tau} \le H^{\tau}$ and
$(\phi^{\tau})_L = \phi(1)\lambda^{\tau} = \phi(1)\lambda$, for any $\tau
\in G$.  It follows that the induced character $\phi^G$ restricts to
$(\phi^G)_L = \phi^G(1)\lambda$ on $L$. So Proposition
\ref{UnderZetaPsiProp} for $\psi = \phi^G$ tells us that $L \le \Z(\phi^G)$
and $\lambda = (\zetaQ{\phi^G})_L$. This, \eqref{KerZetaPsi}, and
\eqref{KerPhiG} imply that  $\Ker(\lambda) \le \Ker(\zetaQ{\phi^G}) =
\Ker(\phi^G) \le \Ker(\phi)$. Thus the present proposition holds.
\end{proof}

Let $H$ be any subgroup of $G$. We say that an irreducible character $\phi
\in \Irr(H)$ \emph{lies under} an irreducible character $\psi \in \Irr(G)$
if $\phi$ is a constituent of the restriction $\psi_H$ of $\psi$ to $H$. By
the Frobenius Reciprocity Law this happens if and only if $\psi$ \emph{lies
over} $\phi$, in the sense that $\psi$ is an irreducible constituent of the
induction $\phi^G$ of $\phi$ to $G$. We indicate that $\phi$ lies under
$\psi$ (or, equivalently,  that $\psi$ lies over $\phi$) by writing $\phi
\le \psi$. We denote by $\Irr(\,G \mid \phi\,)$ the set of all irreducible
characters of $G$ lying over a fixed $\phi \in \Irr(H)$, and by
$\Irr(\,\psi \mid H\,)$ the set of all $\phi \in \Irr(H)$ lying under a
fixed $\psi \in \Irr(G)$. We write $\Lin(\,G \mid \phi\,)$ and 
$\Lin(\,\psi \mid H \,)$ for the subsets
\begin{equation} \mylabel{LinCharSets} \Lin(\, G \mid \phi\,) = \Lin(G)
\cap \Irr(\,G \mid \phi\,) \quad \text{\and} \quad \Lin(\,\psi \mid H\,) =
\Lin(H) \cap \Irr(\,\psi \mid H\,)
\end{equation}
of all linear characters in those two sets. Notice that $\Lin(\,G \mid
\phi\,)$ is empty unless $\phi$ is a linear character, in which case
$\Lin(\,G \mid \phi\,)$ consists of all extensions, if any, of $\phi$ to
linear characters of $G$.

If $K$ is a normal subgroup of $G$, then Clifford theory tells us that
$\Irr(\,\psi \mid K\,)$ is a single $G$-conjugacy class in $\Irr(K)$, for
any fixed $\psi \in \Irr(G)$. It also tells us that induction from
$G(\theta)$ to $G$ is bijection of $\Irr(\,G(\theta) \mid \theta\,)$ onto
$\Irr(\,G \mid \theta\,)$, for any fixed $\theta \in \Irr(K)$. If
$\psi \in \Irr(\,G \mid \theta\,)$, then $\psi_{\theta}$ will denote the
unique character in $\Irr(\,G(\theta) \mid \theta\,)$ from which $\psi =
(\psi_{\theta})^G$ is induced. We call $\psi_{\theta}$ the
\emph{$\theta$-Clifford correspondent} of $\psi$.  The Clifford
correspondent $\psi_{\theta}$ can also be characterized as the only
irreducible character of $G(\theta)$ lying both over $\theta$ and under
$\psi$.  Furthermore, its restriction $(\psi_{\theta})_K$ to $K$ is a
multiple $m\theta$ of $\theta$, for some strictly positive integer $m$.

Suppose that a subgroup $H \le G$ stabilizes a linear character $\kappa$ of
some subgroup $K \le G$. Then $H$ normalizes $K$, so that the product $HK =
KH$ is a subgroup of $G$.  Furthermore, $\kappa$ restricts to an
$H$-invariant linear character $\kappa_{H \cap K}$ of $H \cap K \normaleq
H$.  Clifford theory tells us that any $\phi \in \Irr(\,H \mid \kappa_{H
\cap K}\,)$ restricts to a multiple $\phi(1)\kappa_{H \cap K}$ of
$\kappa_{H \cap K}$ on $H \cap K$. In view of Proposition
\ref{UnderZetaPsiProp}, it follows that any complex matrix
representation $A$ of $H$ affording $\phi$ sends any $\sigma \in H \cap K$
to $\kappa(\sigma)I$. Because $\kappa \in \Lin(K)$ is $H$-invariant, we
conclude that the map
\begin{subequations} \mylabel{HKEqs}
\begin{equation} \mylabel{HKRep} \sigma\tau \mapsto \kappa(\tau)A(\sigma),
\end{equation}
for $\sigma \in H$ and $\tau \in K$, is a well defined irreducible
representation of $HK$. We denote by $\phi * \kappa$ the irreducible
character of $HK$ afforded by this representation. Then
\begin{equation} \mylabel{HKChar}
(\phi * \kappa)(\sigma\tau) = \phi(\sigma)\kappa(\tau) \in \CC,
\end{equation}
\end{subequations}
for any $\sigma \in H$ and $\tau \in K$.  

When $\kappa$ is faithful, the subgroup $K$ is cyclic and central in $HK$, which is the central product $H * K$ of $H$ and $K$ with $H \cap K$ amalgamated. In that case $\phi * \kappa$ is the central product of the two characters $\phi$ and $\kappa$. The case of general $\kappa$ can be reduced to that of faithful $\kappa$ by passing to the factor group $HK/\Ker(\kappa)$.  Like central products, the product $\phi * \kappa$ satisfies
\begin{proposition} \mylabel{HKProp} If\/ $K \le G$ and $H \le G(\kappa)$,
for some $\kappa \in \Lin(K)$, then the map $\phi \mapsto \phi * \kappa$ is
a bijection of\/ $\Irr(\,H \mid \kappa_{H \cap K}\,)$ onto $\Irr(\,HK \mid
\kappa\,)$. The inverse bijection is restriction from $HK$ to $H$.
\end{proposition}
\begin{proof} It is clear from \eqref{HKChar} that $\phi * \kappa$, for any $\phi \in \Irr(\,H \mid \kappa_{H \cap K}\,)$, restricts to $\phi$ on $H$ and to $\phi(1)\kappa$ on $K$. So $\phi * \kappa$ lies in $\Irr(\,HK \mid \kappa\,)$ and determines $\phi$. 

On the other hand, any $\psi \in \Irr(\,HK \mid \kappa\,)$ restricts to a multiple $\psi(1)\kappa$ of the $HK$-invariant linear character $\kappa$ on $K$. In view of \eqref{PsiZetaPsi} and Proposition
\ref{UnderZetaPsiProp}, this implies that $\psi$ restricts to some $\psi_H \in \Irr(\,H \mid \kappa_{H \cap K}\,)$ such that $\psi_H * \kappa = \psi$. Thus the proposition holds.
\end{proof}

\section{ Triples } \mylabel{Triples} 

Recall from the introduction that $\mfrt$ is the family of all triples $T =
(G,N,\psi)$ such that $G$ is a finite group, $N$ is a normal subgroup of
$G$, and $\psi$ is an irreducible character of $N$. We call $G$, $N$, and
$\psi$ the \emph{ambient group}, \emph{normal subgroup}, and
\emph{character}, respectively, in the triple $T$. As in the
introduction, we define the \emph{kernel} $\Ker(T)$ of $T$, the
\emph{center} $\Z(T)$ of $T$, and the \emph{central
character} $\zetaQ{T}$ for $T$ to be the corresponding objects
\begin{subequations} \mylabel{SQObjs}
\begin{align} 
\Ker(T) &= \Ker(\psi^G) = \Core_G(\Ker(\psi)) \le \Ker(\psi) \le N,
\mylabel{KerQ} \\
\Z(T) &= \Z(\psi^G) \le \Z(\psi) \le N, \quad \text{and} \mylabel{ZQ} \\
\zetaQ{T} &= \zetaQ{\psi^G} = (\zetaQ{\psi})_{\Z(T)} \mylabel{ZetaQ} 
\end{align}
\end{subequations}
for the character $\psi^G$ of $G$ induced by $\psi \in \Irr(N)$. So both
$\Ker(T)$ and $\Z(T)$ are normal subgroups of $G$ contained in $N$, while
the linear character $\zetaQ{T}$ of $\Z(T)$ is $G$-invariant with the same
kernel
\begin{equation} \mylabel{KerZetaQ} \Ker(\zetaQ{T}) = \Ker(\zetaQ{\psi^G})
= \Ker(\psi^G) = \Ker(T)
\end{equation}
as $T$.  Hence the factor group $\Z(T)/\Ker(T)$ is cyclic and central in
$G/\Ker(T)$. 

For the rest of this paper we fix a an arbitrary triple $T = (G,N,\psi) \in
\mfrt$, as well as the objects $G$, $N$ and $\psi$ in $T$. A different
characterization of $\Z(T)$ is given in
\begin{proposition} \mylabel{UnderZetaQProp} The center $\Z(T)$ is the
largest normal subgroup $L$ of\/ $G$ contained in $N$ such that $\psi \in
\Irr(N)$ lies over some $G$-invariant linear character $\lambda$ of\/ $L$.
Any such $\lambda$ is the restriction $(\zetaQ{T})_L$ of\/ $\zetaQ{T}$ to
the corresponding $L$. Furthermore, the restriction $\psi_L$ of\/ $\psi$ to
$L$ is a multiple $\psi(1)\lambda$ of\/ $\lambda$. Hence $\Ker(\lambda) \le
\Ker(\zetaQ{T}) = \Ker(T) \le \Ker(\psi)$.
\end{proposition}
\begin{proof} If $\lambda$ is a $G$-invariant linear character of some
normal subgroup $L \normaleq G$ contained in $N$, then Clifford theory for
$L \normaleq N$ implies that $\psi \in \Irr(N)$ lies over $\lambda$ if and
only if $\psi_L = \psi(1)\lambda$. The proposition follows immediately from
this, \eqref{SQObjs}, and Proposition \ref{UnderZetaPhiGProp}.
\end{proof}

We shall often encounter a situation where $L$ and $M$ are normal subgroups of $G$, contained in $N$,
such that $[L,M] \le \Z(T)$. In symbols, this just says that
\begin{subequations} \mylabel{CEqs}
\begin{equation} \mylabel{CLMConds} L,M \normaleq G, \quad L,M \le N,\quad
\text{and} \quad [L,M] \le \Z(T).
\end{equation}
When this happens we may compose commutation with the linear character
$\zetaQ{T}$ of $\Z(T)$ to obtain a function $c = c_T$ from the cartesian
product $L \times M$ to $\CC^{\times}$, sending $(\rho, \sigma) \in L
\times M$ to
\begin{equation} \mylabel{C}  c(\rho,\sigma) =  \zetaQ{T}([\rho,\sigma])
\in \CC^{\times},
\end{equation}
for any $\rho \in L$ and $\sigma \in M$. The commutator identities
$[\pi\rho,\sigma] = [\pi, \sigma]^{\rho}[\rho,\sigma]$ and $[\rho,
\sigma\tau] = [\rho, \tau][\rho,\sigma]^{\tau}$, together with the
$G$-invariance of the linear character $\zetaQ{T}$, imply that $c$ is
\emph{bilinear}, in the sense that
\begin{equation} \mylabel{BilinearC}  c(\pi\rho,\sigma) =
c(\pi,\sigma)c(\rho,\sigma) \in \CC^{\times}  \quad
\text{and} \quad c(\rho, \sigma\tau) = c(\rho,\sigma)c(\rho,\tau) \in
\CC^{\times}, 
\end{equation}
for any $\pi, \rho \in L$ and $\sigma, \tau \in M$. The $G$-invariance of
$\zetaQ{T}$ also implies that $c$ is \emph{$G$-invariant}, in the sense
that
\begin{equation} \mylabel{GInvC}  c(\rho^{\tau}, \sigma^{\tau}) =
c(\rho,\sigma) \in \CC^{\times},
\end{equation}
\end{subequations}
for all $\rho \in L$, $\sigma \in M$ and $\tau \in G$. 

\begin{proposition} \mylabel{CZPerpProp} In the above situation both
$c(\Z(T) \cap L, M)$ and $c(L, \Z(T) \cap M)$ are $1$. Hence $c$ induces a $G$-invariant
bilinear function $\bar c$ from $(L\Z(T)/\Z(T))\times (M\Z(T)/\Z(T))$ to
$\CC^{\times}$, sending the cosets $\rho \Z(T) \in L\Z(T)/\Z(T)$ and
$\sigma \Z(T) \in M\Z(T)/\Z(T)$ to
\begin{equation} \mylabel{CBar} \bar c(\rho\Z(T), \sigma\Z(T)) = c(\rho,
\sigma) = \zetaQ{T}([\rho, \sigma]) \in \CC^{\times}, 
\end{equation}
for any $\rho \in L$ and $\sigma \in M$.
\end{proposition}
\begin{proof} If $\rho \in \Z(T) \cap L$ and $\sigma \in M$, then $[\rho,
\sigma] = \rho^{-1}\rho^{\sigma}$,where both $\rho^{-1}$ and
$\rho^{\sigma}$ lie in $\Z(T) \normaleq G$. Since $\zetaQ{T} \in
\Lin(\Z(T))$ is $G$-invariant, we have
\[  c(\rho, \sigma) = \zetaQ{T}([\rho,\sigma]) =
\zetaQ{T}(\rho^{-1}\rho^{\sigma}) =
\zetaQ{T}(\rho^{-1})\zetaQ{T}(\rho^{\sigma}) =
\zetaQ{T}(\rho)^{-1}\zetaQ{T}(\rho)  =1. \]
So $c(\Z(T) \cap L, M) = 1$. The proof that $c(L, \Z(T) \cap M) = 1$ is similar. Thus the first statement of the proposition holds. The remaining statement follows immediately from the first one and (\ref{CEqs}c,d).
\end{proof}

As in the introduction, we define a \emph{subtriple} of $T$ to be any
triple $T' = (G', N' ,\psi')$, such that
\begin{equation} \mylabel{SubTConds} G' \le G, \quad N' = G' \cap N, \quad
\text{and} \quad \psi' \in \Irr(\,\psi \mid N'\,).
\end{equation}
Of course, any such $T'$ also lies in $\mfrt$. We denote by $\ST(T)$ the
finite family of all subtriples of $T$. We often write $T' \le T$, instead
of $T' \in \ST(T)$, to say that $T'$ is a subtriple of $T$.

For the rest of this section we fix the above subtriple $T' = (G',N',\psi')$ of $T$. Since $N'$ is the intersection $G' \cap N$ of $G' \le G$ with $N \normaleq
G$, there is a natural monomorphism $e^T_{T'}$ of the factor group $G'/N'$
into $G/N$, sending the coset $\sigma N' \in G'/N'$ to the coset
\begin{equation} \mylabel{ETT} e^T_{T'}(\sigma N') = \sigma N \in G/N,
\end{equation}
for any $\sigma \in G'$. We call $e^T_{T'}$ the \emph{natural embedding} of
$G'/N'$ in $G/N$.

It is obvious from the above definition that any subtriple $T'' =
(G'',N'',\psi'')$ of $T'$ is also a subtriple of $T$. Furthermore, the
natural embedding of $G''/N''$ into $G/N$ is the composition
\begin{equation} \mylabel{CompETT} e^T_{T''} = e^T_{T'} \circ e^{T'}_{T''}
\colon G''/N'' \to G/N
\end{equation}
 of that of $G'/N'$ in $G/N$ with that of $G''/N''$ in $G'/N'$.

The center and central character of the subtriple $T' \le T$ are related to those
of $T$ by
\begin{proposition} \mylabel{SubTZZetaProp} The intersection $Z = G' \cap
\Z(T)$ is a normal subgroup of\/ $G'$ contained in $N'$. The restriction
$\zeta = (\zetaQ{T})_Z$ is a $G'$-invariant linear character of\/ $Z$ lying
under $\psi' \in \Irr(N')$.  Hence $Z$ is the intersection $\Z(T') \cap
\Z(T)$, and the restriction $(\zetaQ{T'})_Z$ is also equal to $\zeta$.
\end{proposition}
\begin{proof} Because $\Z(T)$ is a normal subgroup of $G$ contained in $N$,
its intersection with $G' \le G$ is a normal subgroup $Z$ of $G'$ contained
in $N' = G' \cap N$. The $G$-invariant linear character $\zetaQ{T}$ of
$\Z(T)$ restricts to a $G'$-invariant linear character $\zeta =
(\zetaQ{T})_Z$ of $Z$.  We know from Proposition \ref{UnderZetaQProp} that
$\psi \in \Irr(N)$ restricts to a multiple $\psi(1)\zetaQ{T}$ of
$\zetaQ{T}$.  So its restriction $\psi_{N'}$ to $N'$
restricts to the same multiple $\psi(1)\zeta$ of $\zeta$. It follows that
the irreducible constituent $\psi'$ of $\psi_{N'}$ also restricts to a
multiple $\psi'(1)\zeta$ of the linear character $\zeta$. Now all the
hypotheses of Proposition \ref{UnderZetaQProp} hold with $T'$, $Z$, and
$\zeta$ in place of $T$, $L$, and $\lambda$, respectively. That
proposition tells us that $Z \le \Z(T')$, and that  $(\zetaQ{T'})_{Z} =
\zeta$. Thus the present proposition holds.
\end{proof}

Conjugation by any $\tau \in G(\psi)$ leaves invariant both $N \normaleq G$
and $\psi \in \Irr(N)$. So it sends $N' = G' \cap N$ to $(N')^{\tau} =
(G')^{\tau} \cap N$, and $\psi' \in \Irr(\,\psi \mid N'\,)$ to
$(\psi')^{\tau} \in \Irr(\,\psi \mid (N')^{\tau}\,)$. Hence it sends the
subtriple $T' = (G',N',\psi')$ of $T$ to the \emph{conjugate} subtriple
\begin{equation} \mylabel{ConjSubT} (T')^{\tau} = ((G')^{\tau},
(N')^{\tau}, (\psi')^{\tau}) \le T.
\end{equation}
Evidently this conjugation is an action of the group $G(\psi)$ on
the set $\ST(T)$ of all subtriples of $T$. 

Conjugation by $\tau$ sends the character $(\psi')^{G'}$ of $G'$ induced by
$\psi'$ to the character $((\psi')^{\tau})^{(G')^{\tau}}$ of $(G')^{\tau}$
induced by $(\psi')^{\tau}$. It follows that it sends the kernel
$\Ker(T')$, the center $\Z(T')$, and the central character $\zetaQ{T'}$ of
$(\psi')^{G'}$ to the kernel
\begin{subequations} \mylabel{ConjSubTEqs}
\begin{equation} \mylabel{ConjKerT} \Ker((T')^{\tau}) = \Ker(T')^{\tau},
\end{equation}
the center
\begin{equation} \mylabel{ConjZT} \Z((T')^{\tau}) = \Z(T')^{\tau},
\end{equation}
and the central character
\begin{equation} \mylabel{ConjZetaT} \zetaQ{(T')^{\tau}} =
(\zetaQ{T'})^{\tau},
\end{equation}
respectively, of $((\psi')^{\tau})^{(G')^{\tau}}$.  Conjugation by $\tau$
also induces an isomorphism of the factor group $G'/N'$ onto
$(G')^{\tau}/(N')^{\tau}$, sending any coset $\bar{\sigma} \in G'/N'$ to
the coset $\bar{\sigma}^{\tau} \in (G')^{\tau}/(N')^{\tau}$. The natural
embeddings \eqref{ETT} of $G'/N'$ and $(G')^{\tau}/(N')^{\tau}$ in $G/N$
carry this isomorphism into conjugation by $\tau N \in G/N$, in the sense
that
\begin{equation} \mylabel{ConjETT}  e^T_{(T')^{\tau}}(\bar{\sigma}^{\tau})
= e^T_{T'}(\bar{\sigma})^{\tau N} \in G/N,
\end{equation}
\end{subequations}
for any $\bar{\sigma} \in G'/N'$. When $\tau$ lies in the subgroup $N$ of
$G(\psi)$, its image $\tau N = 1_{G/N}$ centralizes $G/N$. In that case
conjugation by $\tau$ carries $e^T_{T'}$ to $e^T_{(T')^{\tau}}$.

\section{ Linear reductions } \mylabel{LinearReductions}

Suppose that $L$ and $\lambda$ satisfy
\begin{equation}\mylabel{LLambdaConds} L \normaleq G, \quad L \le N, \quad
\text{and} \quad \lambda \in \Lin(\,\psi \mid L\,). 
\end{equation}
Then the stabilizer $G(\lambda)$ of $\lambda$ in $G$ is a subgroup of $G$,
the stabilizer $N(\lambda)$ of $\lambda$ in $N$ is the intersection
$G(\lambda) \cap N$, and the $\lambda$-Clifford correspondent
$\psi_{\lambda}$ of $\psi \in \Irr(\, N \mid \lambda\,)$ is an irreducible
character of $N(\lambda)$ lying under $\psi$. Hence
\begin{equation} \mylabel{LinRedQuad} T(\lambda) =  (G(\lambda),
N(\lambda), \psi_{\lambda})
\end{equation}
is a subtriple of $T$. We call $T(\lambda)$ the  \emph{$\lambda$-linear
reduction} of $T$.  Notice that
\begin{equation} \mylabel{LambdaUnderZeta} L \le \Z(T(\lambda)) \quad \text{and}  \quad \lambda = (\zeta^{(T(\lambda))})_L.
\end{equation}
This follows immediately from Proposition \ref{UnderZetaQProp} for $T(\lambda)$, since $L$ is a normal subgroup of $G(\lambda)$ contained in $N(\lambda)$, and $\lambda$ is a $G(\lambda)$-invariant linear character of $L$ lying under $\psi_{\lambda}$.

We say that an arbitrary subtriple $T'$ of $T$ is a
\emph{linear reduction} of $T$ if it is equal to $T(\lambda)$ for
some $L$ and $\lambda$ satisfying \eqref{LLambdaConds}. We denote by
$\LR(T)$ the finite family of all linear reductions of $T$. Note that $T$
lies in $\LR(T)$, since we may choose $L = \Z(T)$ and $\lambda = \zetaQ{T}$
in \eqref{LLambdaConds}, in which case $T(\lambda)$ is just $T$.  

A different way of looking at the conditions \eqref{LLambdaConds} is given
in
\begin{proposition} \mylabel{IsotropicL} If\/ $L \normaleq G$ and $L \le
N$, then $\Lin(\,\psi \mid L\,)$ is non-empty if and only if\/ $[L,L] \le
\Ker(T)$. This happens if and ony if $L\Z(T)/\Z(T)$ is an abelian normal subgroup of\/ $G/\Z(T)$ contained in $N/\Z(T)$.
\end{proposition}
\begin{proof}  The subgroup $[L,L]$, like $L$, is normal in $G$ and
contained in $N$. If $[L,L] \le \Ker(T)$, then $[L,L] \le \Ker(\psi)$ by
\eqref{KerQ}. So any character of $L$ lying under $\psi$ is linear. Hence
$\Lin(\,\psi \mid L\,)$ is the non-empty set $\Irr(\,\psi \mid L\,)$ when
$[L,L] \le \Ker(T)$.

Conversely, suppose we have some character $\lambda \in \Lin(\,\psi \mid
L\,)$. Then $[L,L] \le \Ker(\lambda)$. Because $[L,L]$ is normal in $G$, it
follows that $[L,L] \le \Core_G(\Ker(\lambda)) = \Ker(\lambda^G)$. But 
$\Ker(\lambda^G) \le \Ker(\psi^G) = \Ker(T)$, since $\psi \in \Irr(\,N \mid
\lambda\,)$ is an irreducible constituent of $\lambda^N$.  Hence $[L,L] \le
\Ker(T)$ when $\Lin(\,\psi \mid L \,)$ is not empty, and the first statement of the proposition is
proved.  The remaining statement is an immediate consequence of the first one.
\end{proof}

Assume that $L$ and $\lambda$ satisfy \eqref{LLambdaConds}. Since both $L$
and $\Z(T)$ are normal subgroups of $G$ contained in $N$, so are $L\Z(T)$
and $L \cap \Z(T)$.  The restriction $\lambda_{L \cap \Z(T)} \in \Lin(L
\cap \Z(T))$ lies under $\lambda \le \psi$. Since $\psi_{\Z(T)} =
\psi(1)\zetaQ{T}$ by Proposition \ref{UnderZetaQProp}, this implies that $\lambda_{L \cap \Z(T)} = (\zetaQ{T})_{L \cap \Z(T)}$. The $G$-invariant
character $\zetaQ{T}$ is also $L$-invariant. So  we may apply
\eqref{HKEqs}, with $H = L$, $K = \Z(T)$, $\phi = \lambda$, and $\kappa =
\zetaQ{T}$, to obtain a linear character $\lambda * \zetaQ{T}$ of
$L\Z(T)$ with the value
\begin{equation} \mylabel{LambdaZeta} (\lambda * \zetaQ{T})(\sigma\tau) =
\lambda(\sigma)\zetaQ{T}(\tau) \in \CC^{\times}, 
\end{equation}
for any $\sigma \in L$ and $\tau \in \Z(T)$. Evidently $\lambda *
\zetaQ{T}$ extends $\lambda$. So its domain $L\Z(T)$ is a subgroup of
$G(\lambda)$. In particular, $\Z(T) \le L\Z(T)$ is contained in the ambient
group $G(\lambda)$ of $T(\lambda) \le T$. By Proposition
\ref{SubTZZetaProp} this implies that
\begin{equation} \mylabel{LinRedZZeta} \Z(T) \le \Z(T(\lambda)) \quad
\text{and} \quad \zetaQ{T} = (\zetaQ{T(\lambda)})_{\Z(T)}.
\end{equation}

Another consequence of the existence of $\lambda * \zetaQ{T}$ is
\begin{proposition} \mylabel{ZZetaLLambdaProp} Any linear reduction of\/
$T$ has the form $T(\kappa)$, for some $K$ and $\kappa$ satisfying 
\begin{equation} \mylabel{LLambdaPrimeConds} K \normaleq G, \quad \Z(T)
\le K \le N, \quad \text{and} \quad  \kappa \in
\Lin(\,\psi \mid K\,).
\end{equation}
Furthermore, the restriction $\kappa_{\Z(T)}$ is $\zetaQ{T}$, for any
such $K$ and $\kappa$.
\end{proposition}
\begin{proof}  By definition our linear reduction is $T(\lambda)$ for some
$L$ and $\lambda$ satisfying  \eqref{LLambdaConds}.  As we saw above, the
product $K = L\Z(T)$ is a normal subgroup of $G$ contained in $N$. Because
$\psi \in \Irr(N)$ lies over $\lambda \in \Lin(L)$, and restricts to
$\psi(1)\zetaQ{T}$ on $\Z(T)$, it lies over $\kappa = \lambda * \zetaQ{T}
\in \Lin(K)$.  It follows that $K$ and $\kappa$ satisfy
\eqref{LLambdaPrimeConds}.  In particular, they satisfy the equivalent of
\eqref{LLambdaConds}, so that the linear reduction $T(\kappa)$ of $T$ is
defined.

Since $G$ normalizes $L = \Dom(\lambda)$, and stabilizes $\zetaQ{T}$, it
follows from \eqref{LambdaZeta} that $G(\kappa) = G(\lambda * \zetaQ{T})
= G(\lambda)$. Similarly, $N(\kappa) = N(\lambda)$. So the irreducible
character $\psi_{\kappa}$  of $N(\kappa)$ lying
under $\psi$ and over $\kappa = \lambda * \zetaQ{T}$ must equal the
unique irreducible character $\psi_{\lambda}$ of $N(\lambda)$ lying under
$\psi$ and over $\lambda \le \lambda * \zetaQ{T}$.  Therefore $T(\kappa)
= (G(\kappa), N(\kappa), \psi_{\kappa})$ is equal to $T(\lambda) =
(G(\lambda), N(\lambda), \psi_{\lambda})$, and the first statement in the
proposition is proved.

Now let $K$ and $\kappa$ be any objects satisfying
\eqref{LLambdaPrimeConds}.  The subgroup $K$ of $G(\kappa)$ is contained
in $N(\kappa)$, and its $G(\kappa)$-invariant linear character
$\kappa$ lies under $\psi_{\kappa}$. So Proposition
\ref{UnderZetaQProp}, for the triple $T(\kappa)$, tells us that $K \le
\Z(T(\kappa))$ and $\kappa = (\zetaQ{T(\kappa)})_{K}$.   Since
$\Z(T) \le K$ by hypothesis, this and \eqref{LinRedZZeta}, with $\kappa$
in place of $\lambda$, imply that
\[  \kappa_{\Z(T)} = ((\zetaQ{T(\kappa)})_{K})_{\Z(T)} =
(\zetaQ{T(\kappa)})_{\Z(T)} = \zetaQ{T}. \]
Hence the remaining statement in the proposition holds.
\end{proof}

We say that a linear reduction of $T$ is \emph{proper} if it is not
equal to $T$.  We say that $T$ is \emph{linearly reducible} if it has some
proper linear reduction, and \emph{linearly irreducible} otherwise.
\begin{proposition} \mylabel{PropLinRedProp} If\/ $L$ and $\lambda$ satisfy
\eqref{LLambdaConds}, then $T(\lambda)$ is a proper linear reduction of\/
$T$ if and only if\/ $G(\lambda) < G$. This happens if and only if\/
$\Z(T)$ does not contain $L$. Thus $T$ is linearly irreducible if and only
if\/ $\Z(T)$ contains every normal subgroup $L$ of\/ $G$ satisfying both $L
\le N$ and $[L,L] \le \Ker(T)$.
\end{proposition}
\begin{proof} Suppose that $G(\lambda) = G$. Then $N(\lambda) = G(\lambda)
\cap N = N$, and $\psi_{\lambda}$ must be the unique character $\psi \in
\Irr(\,\psi \mid N\,)$.  Hence $T(\lambda) = T$ in this case. Since
$T(\lambda) \ne T$ when $G(\lambda) < G$, this proves the first statement
of the proposition.

By definition $L$ is a normal subgroup of $G$ contained in $N$, and
$\lambda$ is a linear character of $L$ lying under $\psi \in \Irr(N)$.
Proposition \ref{UnderZetaQProp} implies that $G$ stabilizes $\lambda$ if
and only if $L$ is a subgroup of $\Z(T)$. The second statement of the
proposition follows from this and the first statement.

The triple $L$ is linearly irreducible if and only if it is equal to
$T(\lambda)$ whenever $L$ and $\lambda$ satisfy \eqref{LLambdaConds}. By
the second statement of the proposition this happens if and only if $\Z(T)$
contains $L$, for all such $L$ and $\lambda$. In view of Proposition
\ref{IsotropicL} this is equivalent to the final statement of the 
proposition.
\end{proof}

We can restate the above proposition in terms of factor groups modulo $\Ker(T)$.
\begin{proposition} \mylabel{AltLinIrrProp} The arbitrary triple $T = (G,N,\psi) \in 
\mfrt$ is linearly irreducible if and only if the cylic central subgroup
$\Z(T)/\Ker(T)$ of\/ $G/\Ker(T)$ is the largest abelian normal subgroup
of\/ $G/\Ker(T)$ contained in $N/\Ker(T)$.  In that case $\Z(T)/\Ker(T)$ is
the center $\Z(N/\Ker(T))$ of\/ $N/\Ker(T)$.
\end{proposition}
\begin{proof}  As we noted after \eqref{KerZetaQ}, the subgroup
$\Z(T)/\Ker(T)$ is always cyclic and central in $G/\Ker(T)$. It is also
contained in $N/\Ker(T)$ by \eqref{ZQ}.  The normal subgroups $\bar L$ of
$G/\Ker(T)$ contained in $N/\Ker(T)$ are just the images $\bar L =
L\Ker(T)/\Ker(T)$ of the normal subgroups $L$ of $G$ contained in $N$.  
Clearly $\bar L$ is abelian if and only if any corresponding $L$ satisfies
$[L,L] \le \Ker(T)$. The first statement of the proposition follows from
these remarks and Proposition \ref{PropLinRedProp}. The remaining statement
is an immediate consequence of the first one.
\end{proof}

An elementary exercise in Clifford theory gives us
\begin{proposition} \mylabel{CompoundLinRedProp} If\/ $L$ and $\lambda$
satisfy \eqref{LLambdaConds}, and if $\kappa$ is the restriction
$\lambda_K$ of\/ $\lambda \in \Lin(L)$ to some normal subgroup $K$ of\/ $G$
contained in $L$, then the $\kappa$-linear reduction $T(\kappa)$ of\/ $T$
is defined. Furthermore, the $\lambda$-linear reduction $T(\kappa, \lambda)
= [T(\kappa)](\lambda)$ of\/ $T(\kappa)$ is defined and equal to
$T(\lambda)$.
\end{proposition} 
\begin{proof} The normal subgroup $K$ of $G$ is contained in $N$, since it
is contained in $L \le N$. Its linear character $\kappa = \lambda_K$ lies
under $\psi \in \Irr(N)$, because $\lambda$ does. Hence
\eqref{LLambdaConds} holds with $K$ and $\kappa$ in place of $L$ and
$\lambda$, respectively. So the $\kappa$-linear reduction $T(\kappa) =
(G(\kappa), N(\kappa), \psi_{\kappa})$ is defined.

Clearly $L$ fixes the restriction $\kappa$ of its character $\lambda$ to
its normal subgroup $K$. Hence $L$ is a normal subgroup of $G(\kappa)$
contained in $N(\kappa)$. Because $\lambda$ extends $\kappa$, and lies
under $\psi$, it lies under the unique irreducible constituent
$\psi_{\kappa}$ of $\psi_{\N(\kappa)}$ lying over $\kappa$. Thus
\eqref{LLambdaConds} holds with $T(\kappa)$ in place of $T$. Hence the
$\lambda$-linear reduction $T(\kappa, \lambda) = [T(\kappa)](\lambda) =
(G(\kappa,\lambda), N(\kappa,\lambda), (\psi_{\kappa})_{\lambda})$
is defined.

Since $K$ is a normal subgroup of $G$, any element of $G$ fixing $\lambda$
must fix the restriction $\kappa = \lambda_K$.  Hence $G(\kappa, \lambda) =
G(\lambda)$. It follows that $N(\kappa, \lambda) = G(\kappa, \lambda) \cap
N = G(\lambda) \cap N = N(\lambda)$. So $(\psi_{\kappa})_{\lambda}$ is an
irreducible character of $N(\lambda)$ lying over $\lambda$ and under
$\psi_{\kappa} \le \psi$. Thus it must be the unique irreducible character
$\psi_{\lambda}$ of $N(\lambda)$ lying over $\lambda$ and under $\psi$.
Therefore $T(\kappa, \lambda)$ is $T(\lambda) = (G(\lambda), N(\lambda),
\psi_{\lambda})$, and the  proposition is proved.
\end{proof}

\section{ Multilinear reductions } \mylabel{MultLinReductions}

As in the introduction, we say that that a subtriple $T' \le T$ is a
\emph{multilinear reduction} of $T$ if there is some finite chain $T_0,
T_1, \dots, T_n$ of subtriples $T_i \le T$, starting with $T_0 = T$ and
ending with $T_n = T'$, such that each $T_i$, for $i = 1,2,\dots,n$, is
a linear reduction of its predecessor $T_{i-1}$.  In that case we call
$T_0, T_1,\dots, T_n$ a \emph{linear reducing chain} from $T$ to $T'$.
We write $\MLR(T)$ for the finite family of all multilinear reductions
of $T$. Of course, the set $\LR(T)$ of all linear reductions of $T$ is a
subset of $\MLR(T)$. In particular, $T \in \LR(T)$ lies in $\MLR(T)$.
It follows immediately from \eqref{LinRedZZeta} that
\begin{equation} \mylabel{MLRZZeta} \Z(T) \le \Z(T') \quad \text{and} \quad
\zetaQ{T} = (\zetaQ{T'})_{\Z(T)}
\end{equation}
for any $T' \in \MLR(T)$.

Another concept defined in the introduction is that of a \emph{linear
limit} of $T$. This is any multilinear reduction $T'$ of $T$ such that
$T'$ is linearly irreducible. We denote by $\LL(T)$ the family of all
linear limits of $T$.  Evidently $\LL(T)$ is a non-empty subset of
$\MLR(T)$. Furthermore, $\LL(T')$ is a non-empty subset of $\LL(T)$, for
any $T' \in \MLR(T)$.

Let $T' =(G',N',\psi')$ be an arbitrary subtriple of $T =
(G,N,\psi)$.  Any linear reduction of $T'$ has the form $T'(\lambda')$, for
some normal subgroup $L'$ of $G'$ contained in $N'$, and some linear
character $\lambda'$ of $L'$ lying under $\psi' \in \Irr(N')$.
Conjugation by any $\tau \in G(\psi)$ sends $L'$ to some normal subgroup
$(L')^{\tau}$ of $(G')^{\tau}$ contained in $(N')^{\tau}$. It also sends
$\lambda'$ to some linear character $(\lambda')^{\tau}$ of $(L')^{\tau}$
lying under $(\psi')^{\tau} \in \Irr((N')^{\tau})$. Since
$(G')^{\tau}((\lambda')^{\tau}) = G'(\lambda')^{\tau}$ and
$((\psi')^{\tau})_{(\lambda')^{\tau}} = ((\psi')_{\lambda'})^{\tau}$, it
follows from \eqref{ConjSubT} and \eqref{LinRedQuad} that
$(T')^{\tau}((\lambda')^{\tau})$ is the conjugate subtriple
\begin{equation} \mylabel{ConjSQLambda} (T')^{\tau}((\lambda')^{\tau}) =
T'(\lambda')^{\tau} 
\end{equation}
of $T$.

Since the conjugate $T^{\tau} = (G^{\tau}, N^{\tau}, \psi^{\tau})$ is equal
to $T = (G,N,\psi)$, for any $\tau \in G(\psi)$, the above equation, together with the definitions of linear reductions, multilinear reductions, and linear limits, implies immediately that
\begin{proposition} \mylabel{ConjLLProp} If\/ $T' \in \ST(T)$ and $\tau \in
G(\psi)$, then conjugation by $\tau$ sends the subsets $\LR(T')$,
$\MLR(T')$ and $\LL(T')$ of\/ $\ST(T)$ onto 
\[ \LR(T')^{\tau} = \LR((T')^{\tau}), \quad \MLR(T')^{\tau} =
\MLR((T')^{\tau})\quad \text{and}\quad \LL(T')^{\tau} = \LL((T')^{\tau}),
\]
respectively. In particular, $\LR(T)$, $\MLR(T)$ and $\LL(T)$ are
$G(\psi)$-invariant subsets of\/ $\ST(T)$.
\end{proposition}

A little exercise in Clifford theory will give us
\begin{proposition} \mylabel{MultLinRedStabProp} If\/ $T' = (G',N',\psi')$ is a multilinear reduction of\/ $T$ \textup{(}in particular, if\/ $T'$ is a linear limit of\/
$T$\textup{)}, then $G'(\psi')N = G(\psi)$. Hence the natural embedding $e^T_{T'}$
of $G'/N'$ into $G/N$ sends $G'(\psi')/N' \le G'/N'$ isomorphically onto
$G(\psi)/N \le G/N$.
\end{proposition}
\begin{proof} Notice that the two statements in the proposition are
equivalent to each other, since the monomorphism $e^T_{T'}$ of $G'/N'$ into
$G/N$ in \eqref{ETT} sends $G(\psi')/N'$ onto $G(\psi')N/N$. Furthermore,
the definition of multilinear reductions as repeated linear
reductions, together with the transitivity \eqref{CompETT} of natural
embeddings, implies that the proposition will hold in general if it holds
in the special case where $T'$ is a linear reduction of $T$. So all we have
to prove is the equality $G'(\psi')N = G(\psi)$ in that case.

Since $T'$ is a linear reduction of $T$, there exist some $L$ and $\lambda$
satisfying \eqref{LLambdaConds} such that $T' = T(\lambda)$. Then  $G'$ is $G(\lambda)$, and $\psi'$ is
$\psi_{\lambda}$. It follows that
\[ G'(\psi') = G'(\psi_{\lambda}) = G(\lambda, \psi_{\lambda}).
\]

By definition the Clifford correspondent $\psi_{\lambda}$ induces the
character $\psi$ of $N$. Since the subgroup $G(\lambda, \psi_{\lambda})$ of
$G$ fixes $\psi_{\lambda}$ and normalizes $N \normaleq G$, it fixes $\psi =
(\psi_{\lambda})^N$. Hence $G(\lambda, \psi_{\lambda}) \le G(\lambda,
\psi)$.  The opposite inclusion comes from the observation that $G(\lambda,
\psi)$ fixes both $\lambda$ and $\psi$, and therefore fixes the unique
$\lambda$-Clifford correspondent $\psi_{\lambda}$ of $\psi$. Thus we have
 \[  G'(\psi') = G(\lambda, \psi_{\lambda}) = G(\lambda, \psi). \]

Since $L$ is a normal subgroup of $G$, conjugation by any element $\tau \in
G(\psi)$ leaves invariant the set $\Irr(\,\psi \mid L\,)$ of all
irreducible characters of $L$ lying under $\psi \in \Irr(N)$. Clifford
theory for $L \normaleq N$ tells us that $\Irr(\,\psi \mid L\,)$ is a
single $N$-conjugacy class. Because $\lambda$ lies in $\Irr(\,\psi \mid
L\,)$ by \eqref{LLambdaConds}, it follows that $G(\psi)$ is the product
$G'(\psi')N$ of the stabilizer $G'(\psi') = G(\lambda, \psi)$
of $\lambda$ in $G(\psi)$ with $N$. As we saw above, this is enough to
prove the proposition.
\end{proof} 

Any multilinear reduction induces several correspondences of characters.
\begin{proposition} \mylabel{MLRCharCorrProp} If\/ $T' = (G',N',\psi')$ is a multilinear reduction of\/ $T$, and if\/ $H$ is a subgroup of $G$ containing $N$, then $H' = G' \cap H$ is a subgroup of\/ $G'$ containing $N'$. Induction from $H'$ to $H$
is a bijection of\/ $\Irr(\,H' \mid \psi'\,)$ onto $\Irr(\,H \mid \psi\,)$.  
The inverse bijection sends any $\theta \in \Irr(\,H \mid \psi\,)$ to the
unique character $\theta_{T'}$ in $\Irr(\,H' \mid \zeta^{(T')}\,)$ lying under $\theta$.
\end{proposition}
\begin{proof}  The subtriple $T'$ of $T$ satisfies $N' = G'\cap N$ by
\eqref{SubTConds}. The first statement of the proposition is an immediate
consequence of this.  

We first prove the remaining statements of the proposition when $T'$ is a linear reduction of $T$. In this case there exist some $L$ and $\lambda$ satisfying \eqref{LLambdaConds} such that $T' = T(\lambda) = (G(\lambda), N(\lambda), \psi_{\lambda})$.  It follows that $H' =  G(\lambda) \cap H = H(\lambda)$. So Clifford theory for $L \normaleq H$ and $\lambda \in \Irr(L)$ tells us that induction is a bijection of $\Irr(\,H' \mid \lambda\,)$ onto $\Irr(\,H \mid \lambda\,)$, and that the inverse bijection sends any $\theta \in \Irr(\,H \mid \lambda\,)$ to the unique $\theta_{\lambda} \in \Irr(\,H' \mid \lambda\,)$ lying under $\theta$.  Since $L \le N' \le H'$, while $\psi' = \psi_{\lambda} \in \Irr(N')$ lies over $\lambda$ and induces $\psi \in \Irr(N)$, this bijection sends the subset $\Irr(\,H' \mid \psi'\,)$ of $\Irr(\,H' \mid \lambda\,)$ into $\Irr(\,H \mid \psi\,) \subseteq \Irr(\,H \mid \lambda\,)$.

We know from \eqref{LambdaUnderZeta} for $T' = T(\lambda)$ that $\lambda$ is the restriction of $\zeta^{(T')} \in \Lin(\Z(T'))$ to $L \le \Z(T')$. So any character $\theta' \in \Irr(\,H' \mid \zeta^{(T')}\,)$ lying under $\theta$ lies over $\lambda$, and thus must be $\theta_{\lambda}$ by Clifford theory.  Because $\theta$ lies over $\psi \in \Irr(\,N \mid \lambda\,)$, its Clifford correspondent $\theta' = \theta_{\lambda}$ lies over $\psi' = \psi_{\lambda}$. Hence $\theta'$ must be the unique element of $\Irr(\,H' \mid \psi'\,)$ inducing $\theta$. Thus the remaining two statements of the proposition hold when $T'$ is a linear reduction of $T$.

A general $T'$ is the last member $T' = T_n$ in some linear reducing chain $T_0 = T, T_1, \dots, T_n$.  Let $T_i = (G_i,N_i,\psi_i)$ and $H_i = G_i \cap H$, for $i = 0,1,\dots,n$.  Since $T_i$ is a linear reduction of $T_{i-1}$, for each $i = 1,2,\dots,n$, the above arguments tell us that induction is a bijection of $\Irr(\,H_i \mid \psi_i\,)$ onto $\Irr(\,H_{i-1}\mid \psi_{i-1}\,)$, and that the inverse bijection sends any $\theta_{i-1} \in \Irr(\,H_{i-1}\mid \psi_{i-1}\,)$ to the unique $\theta_i \in \Irr(\,H_i \mid \zeta^{(T_i)}\,)$ lying under $\theta_{i-1}$.  Composing these inductive bijections we see that induction is a bijection of $\Irr(\,H' \mid \psi'\,) = \Irr(\,H_n \mid \psi_n\,)$ onto $\Irr(\,H \mid \psi\,) = \Irr(\,H_0 \mid \psi_0\,)$. Furthermore, the inverse bijection sends any $\theta \in \Irr(\,H \mid \psi\,)$ to some $\theta_{T'} \in \Irr(\,H'\mid \zeta^{(T')}\,)$ lying under $\theta$.  

To complete the proof of the proposition we must show that any  $\theta' \in \Irr(\,H' \mid \zeta^{(T')}\,)$ lying under $\theta$ belongs to $\Irr(\,H' \mid \psi'\,)$ and induces $\theta$.  Since $H' = H_n \le H_{n-1} \le \dots \le H_0 = H$, we may choose characters $\theta_i \in \Irr(H_i)$, for $i = 0,1,\dots,n$, such that $\theta' = \theta_n \le \theta_{n-1} \le \dots \le \theta_0 = \theta$.  It follows from \eqref{LinRedZZeta} that $\zeta^{(T)} = \zeta^{(T_0)} \le \zeta^{(T_1)} \le \dots \le \zeta^{(T_n)} = \zeta^{(T')} \le \theta'$. Hence $\theta_i$ lies in $\Irr(\,H_i \mid \zeta^{(T_i)}\,)$ and under $\theta_{i-1}$, for each $i = 1,2,\dots,n$.  By induction on $i$ this implies that $\theta_i$ lies in $\Irr(\,H_i \mid \psi_i\,)$ and induces $\theta_{i-1}$, for each $i = 1,2,\dots,n$. Hence $\theta' = \theta_n$ lies in $\Irr(\,H' \mid \psi'\,) = \Irr(\,H_n \mid \psi_n\,)$ and induces $\theta = \theta_0$. As we noted above, this completes the proof of the proposition.
\end{proof}

\noindent
We call the character $\theta_{T'} \in \Irr(\,H'
\mid \psi'\,)$ defined in the preceding proposition the
\emph{$T'$-correspondent} of a given character $\theta \in \Irr(\,H \mid
\psi\,)$.

\medskip
Suppose that $T'$  is a \emph{proper} multilinear reduction of $T$, i.~e., a multilinear reduction not equal to $T$. Then we may fix some $L$ and $\lambda$ satisfying \eqref{LLambdaConds} such that $T(\lambda)$ is a proper linear reduction of $T$, and $T'$ is a multilinear reduction of $T(\lambda)$. In symbols, this says that
\begin{equation} \mylabel{TLambdaConds} T(\lambda) < T \quad \text{and} \quad T' \in \MLR(T(\lambda)).
\end{equation}
We know from \eqref{LambdaUnderZeta} that $L \le \Z(T(\lambda))$ and $\lambda = (\zetaQ{T(\lambda)})_{L}$. Since $T'$ is a multilinear reduction of $T(\lambda)$, we know from \eqref{MLRZZeta} that $\Z(T(\lambda)) \le \Z(T')$ and $\zetaQ{T(\lambda)} = (\zetaQ{T'})_{\Z(T(\lambda))}$. Hence
\begin{equation} \mylabel{LambdaUnderZetaPrime} L \le \Z(T') \quad \text{and} \quad \lambda = (\zetaQ{T'})_{L}. \end{equation}

The multilinear reduction $T'$ is actually determined by its central character $\zetaQ{T'}$.
\begin{proposition} \mylabel{ZetaDetProp} The ambient group $G'$ in any multilinear reduction $T' = (G',N',\psi')$ of\/ $T= (G,N,\psi)$ is the stabilizer $G(\zetaQ{T'})$ of the central character $\zetaQ{T'}$ for $T'$. The normal subgroup $N'$ in $T'$ is $N(\zetaQ{T'})$. The character $\psi'$ in $T'$ is the unique irreducible character of\/ $N'$ lying both under $\psi \in \Irr(N)$ and over $\zetaQ{T'}$. In this way $T$ and $\zetaQ{T'}$ determine $T'$ completely.
\end{proposition}
\begin{proof} We prove this by induction on the finite index $[G:G']$ of $G'$ in $G$. If this index is $1$, then $G' = G$ and $T' = T$. In this case the proposition is trivial, since $\psi' = \psi$ lies over $\zetaQ{T'}  =\zetaQ{T}$. So we may assume that $[G:G'] > 1$, and that the proposition holds for all strictly smaller values of this index.

Now we may fix $L$ and $\lambda$ satisfying \eqref{LLambdaConds} such that \eqref{TLambdaConds} holds. The stabilizer $G(\zetaQ{T'})$ of $\zetaQ{T'}$ in $G$ normalizes $L \normaleq G$ and fixes $\lambda = (\zetaQ{T'})_{L}$. Hence it is contained in $G(\lambda)$. So $G(\zetaQ{T'})$ is equal to $[G(\lambda)](\zetaQ{T'})$. Since $T'$ is a multilinear reduction of $T(\lambda)$, and $G(\lambda) < G$, the induction hypothesis tells us that $[G(\lambda)](\zetaQ{T'}) = G'$. Thus $G'$ is $G(\zetaQ{T'})$. It follows that $N' = G' \cap N$ is $N(\zetaQ{T'})$.

We may choose the subgroup $H$ in Proposition \ref{MLRCharCorrProp} to be $N$. Then $H' = G' \cap H = N'$. Furthermore $\psi$ and $\psi'$ are the only characters in $\Irr(\,N \mid \psi\,)$ and $\Irr(\,N' \mid \psi'\,)$, respectively. So that proposition tells us that $\psi'$ is the unique character in $\Irr(N')$ lying both over $\zetaQ{T'}$ and under $\psi$. This completes the proof of the proposition.
\end{proof}

We are not really going to need the rest of the results in this section, but they are too remarkable to omit.  We start with
\begin{lemma} \mylabel{RedLemma} Suppose that $T' = (G', N', \psi')$ is a subtriple of\/ $T = (G,N,\psi)$, that $L$ is a normal subgroup of\/ $G$ contained in $\Z(T')$, and that $\lambda$ is the restriction of\/ $\zetaQ{T'}$ to $L$. Then the linear reduction $T(\lambda)$ is defined, and has $T'$ as a subtriple.
\end{lemma}
\begin{proof} The normal subgroup $N'$ in $T' \le T$ is contained in the normal subgroup $N$ in $T$. So $L \le \Z(T') \le N' \le N$. The character $\psi'$ in $T' \le T$ lies under the character $\psi$ in $T$. Hence $\lambda \le \zetaQ{T'} \le \psi' \le \psi$.  Thus $L$ and $\lambda$ satisfy the conditions \eqref{LLambdaConds}, and the linear reduction $T(\lambda) = (G(\lambda), N(\lambda), \psi_{\lambda})$ is defined.

The ambient group $G'$ in $T'$ fixes the central character $\zetaQ{T'}$ of that subtriple. Furthermore, $G' \le G$ normalizes $L \normaleq G$. So $G'$ fixes the restriction $\lambda$ of $\zetaQ{T'}$ to $L$. Hence $G' \le G(\lambda)$. It follows that $N' = G' \cap N = G' \cap G(\lambda) \cap N = G' \cap N(\lambda)$. Thus $N' \le N(\lambda) \le N$. Since the character $\psi' \in \Irr(N')$ lies under $\psi \in \Irr(N)$, there is some character $\phi \in \Irr(N(\lambda))$ such that $\psi' \le \phi \le \psi$. Then $\phi$ lies over $\lambda \le \psi'$, and so must be the $\lambda$-Clifford correspondent $\psi_{\lambda}$ of $\psi$. Therefore $T'$ is a subtriple of $T(\lambda)$, and the lemma is proved.
\end{proof}

One striking consequence of the above lemma is
\begin{proposition} \mylabel{MLRSubProp}  If\/ $T'$ is a multilinear reduction of\/ $T$, and is also a subtriple of some subtriple $\sq T \le T$, then $T'$ is a multilinear reduction of\/ $\sq T$.
\end{proposition}
\begin{proof}  As usual, we set $T = (G,N,\psi)$, $T' = (G',N',\psi')$ and $\sq T = (\sq G, \sq N,\sq{\psi})$. We shall prove the proposition by  induction on the finite index $[G:G']$. If $[G:G'] = 1$, then $G' = G$ and $T' = T$. Since $T' \le \sq T \le T$, this forces $\sq T$ to equal $T$. So the proposition holds trivially when $[G:G'] = 1$.

From now on we assume that $[G:G'] > 1$, and that the proposition holds for all strictly smaller values of $[G:G']$.  We may choose $L$ and $\lambda$ satisfying \eqref{LLambdaConds} such that \eqref{TLambdaConds} holds. We know from \eqref{LambdaUnderZetaPrime} that $\lambda$ is the restriction of $\zetaQ{T'} \in \Lin(\Z(T'))$ to $L \le \Z(T')$. Since $T' \le \sq T$, we have $L \le N' \le \sq N$ and $\lambda \le \psi' \le \sq{\psi}$. Furthermore, $L \normaleq G$ is normal in $\sq G \le G$. Thus the linear reduction $\sq T(\lambda) = (\sq G(\lambda), \sq N(\lambda), \sq{\psi}_{\lambda})$ is defined. Lemma \ref{RedLemma}, with $\sq T$ in place of $T$, tells us that $T'$ is a subtriple of $\sq T(\lambda)$.  Evidently $\sq G(\lambda)$ is a subgroup of $G(\lambda)$. Hence  $\sq N(\lambda)$ is a subgroup of $N(\lambda)$.  The $\lambda$-Clifford correspondent $\sq{\psi}_{\lambda}$ of $\sq{\psi} \le \psi$ must lie under the $\lambda$-Clifford correspondent $\psi_{\lambda}$ of $\psi$. Therefore $\sq T(\lambda)$ is a subtriple of $T(\lambda)$. Now all the hypotheses of the proposition are satisfied with $\sq T(\lambda)$ and $T(\lambda)$ in place of $\sq T$ and $T$, respectively. Since $[G(\lambda) : G'] < [G : G']$, we know by induction that $T'$ is a multilinear reduction of $\sq T(\lambda)$. Therefore $T'$ is a multilinear reduction of $\sq T$, and the proposition is proved.
\end{proof}

When $T'$ is a proper multilinear reduction of $T$ there is a considerable choice of $L$ and $\lambda$ satisfying  \eqref{LLambdaConds} and \eqref{TLambdaConds}. Surprisingly, a canonical choice is possible.
\begin{proposition} \mylabel{CanonLProp} Suppose that a multilinear reduction $T'$ of\/ $T = (G,N,\psi)$ is not equal to $T$. Let $K$ be the $G$-core of\/ $\Z(T')$, and $\kappa$ be the restriction to $K$ of\/ $\zetaQ{T'} \in \Lin(\Z(T'))$. Then the linear reduction $T(\kappa)$ of\/ $T$ is defined and proper, and has $T'$ as a multilinear reduction.
\end{proposition}
\begin{proof} Set $T' = (G',N',\psi')$.  Clearly $K = \Core_G(\Z(T'))$ is a normal subgroup of $G$ contained in $\Z(T')$. So Lemma \ref{RedLemma} tells us that $T(\kappa)$ is defined and has $T'$ as a subtriple.  Proposition \ref{MLRSubProp} for $\sq T = T(\kappa)$ now says that $T'$ is a multilinear reduction of $T(\kappa)$. So the only problem is to show that $T(\kappa)$ is a proper linear reduction of $T$, i.e., that $G(\kappa) < G$.

Since $T'$ is not equal to $T$, there are some $L$ and $\lambda$ satisfying \eqref{LLambdaConds} such that \eqref{TLambdaConds} holds. Then $L \le \Z(T')$ by \eqref{LambdaUnderZetaPrime}. Since $L$ is a normal subgroup of $G$, this implies that $L \le K = \Core_G(\Z(T'))$.
The restriction $\lambda$ of $\zetaQ{T'}$ in \eqref{LambdaUnderZetaPrime} is also a restriction of $\kappa = (\zetaQ{T'})_K$. Because $L$ is normal in $G$, this implies that $G(\kappa) \le G(\lambda) < G$. So $T(\kappa) < T$, and the proposition is proved.
\end{proof}

Obviously the above proposition can be compounded to obtain a canonical linear reducing chain from $T$ to any multilinear reduction $T'$ of $T$.  If $T' = T$, then this chain consists only of $T_0 = T$. Otherwise we let $T_1$ be the canonical linear reduction $T(\kappa)$ in the proposition. Then $T'$ is a multilinear reduction of $T_1$. If $T' = T_1$ we stop. Otherwise we take $T_2$ to be the canonical linear reduction obtained by applying the proposition with $T_1$ in place of $T$. We continue in this manner until we arrive at $T'$. This must happen eventually, because $T' \le \dots < T_2 < T_1 < T_0 = T$ by construction, and our groups are all finite.

\section{ Covers } \mylabel{Covers} 

As in the introduction, we say that a subtriple $\sq T = (\sq G, \sq N, \sq{\psi})$ of our arbitrary triple $T = (G,N,\psi) \in \mfrt$ \emph{covers $T$ modulo $\Z(T)$} if the subgroup $\sq G \le G$ covers the factor
group $G/\Z(T)$, in the usual sense that
\begin{equation} \mylabel{CoverG} G = \sq G\Z(T).
\end{equation}
In that case the subgroup $\sq G$ completely determines the
subtriple $\sq T$.
\begin{proposition} \mylabel{CoverTProp} There is a one to one
correspondence between all subgroups $\sq G \le G$ covering $G/\Z(T)$ and
all subtriples $\sq T \le T$ covering $T$ modulo $\Z(T)$. Here two such $\sq G$
and $\sq T$ correspond if and only if
\begin{equation} \mylabel{CoverT} \sq T = (\sq G, \sq G \cap N, \psi_{\sq G
\cap N}).
\end{equation}
\end{proposition}
\begin{proof} Suppose that $\sq G$ is a subgroup of $G$ covering $G/\Z(T)$.
Since $\Z(T) \le N \normaleq G$, it follows from \eqref{CoverG} that $\sq G
\cap N $ is a normal subgroup of $\sq G$ such that $N = (\sq G \cap
N)\Z(T)$. Because $\zetaQ{T} \in \Lin(\Z(T))$ is $G$-invariant, we may
apply Proposition \ref{HKProp} with $H = \sq G \cap N$, $K = \Z(T)$, and
$\kappa = \zetaQ{T}$ to conclude that $\psi \in \Irr(\,(\sq G \cap N)\Z(T)
\mid \zetaQ{T}\,)$ restricts to an irreducible character $\psi_{\sq G \cap
N}$ of $\sq G \cap N$. Hence the right side of \eqref{CoverT} is a
subtriple covering $T$ modulo $\Z(T)$.

Let $\sq T = (\sq G, \sq N, \sq{\psi})$ be any subtriple covering $T$
modulo $\Z(T)$. By the definition of subtriples $\sq G$ is a subgroup of
$G$, while $\sq N$ is the intersection $\sq G \cap N$, and $\sq{\psi} \in
\Irr(\sq N)$ lies under $\psi \in \Irr(N)$.  Furthermore, $\sq G$ covers
$G/\Z(T)$. As we saw in the preceding paragraph, this implies that
$\psi_{\sq N}$ is an irreducible character of $\sq N$. So $\psi_{\sq N}$
must equal its irreducible constituent $\sq{\psi}$. Therefore $\sq T$
satisfies \eqref{CoverT}, and the proposition is proved.
\end{proof}

For the rest of this section we fix a subtriple $\sq T = (\sq G, \sq N,
\sq{\psi})$ covering $T$ modulo $\Z(T)$, as well as the objects $\sq G \le
G$, $\sq N = \sq G \cap N$ and $\sq{\psi} = \psi_{\sq N}$ in this
subtriple. Since $\Z(T) \le N \normaleq G$, it follows from \eqref{CoverG}
that $G = \sq GN$. Hence the natural embedding $e^T_{\sq T}$ of $\sq
G/\sq N$ in $G/N$ is an isomorphism
\begin{subequations} \mylabel{CoverEqs}
\begin{equation} \mylabel{CoverETT} e^T_{\sq T} \colon \sq G/\sq N
\overset{\sim}{\longrightarrow} \sq GN/N = G/N
\end{equation}
of groups. Since $\sq{\psi} \in \Irr(\sq N)$ is
the restriction $\psi_{\sq N}$ of $\psi \in \Irr(N)$, this implies  that the character $\psi^G$ of $G$ induced by
$\psi$ restricts to the character
\begin{equation} \mylabel{CoverPsiG}  \sq{\psi}^{\sq G} = (\psi^G)_{\sq G}
\end{equation}
of $\sq G$ induced by $\sq{\psi}$. It follows that the kernel
$\Ker(\sq T)$ of $\sq{\psi}^{\sq G}$ is the intersection
\begin{equation} \mylabel{CoverKer} \Ker(\sq T) = \sq G \cap \Ker(T)
\end{equation}
of $\sq G$ with the kernel $\Ker(T)$ of $\psi^G$, that the center $\Z(\sq
T)$ of $\sq{\psi}^{\sq G}$ is the intersection
\begin{equation} \mylabel{CoverZ} \Z(\sq T) = \sq G \cap \Z(T)
\end{equation}
of $\sq G$ with the center $\Z(T)$ of $\psi^G$, and that the central
character $\zetaQ{\sq T}$ of $\sq{\psi}^{\sq G}$ is the restriction
\begin{equation} \mylabel{CoverZeta} \zetaQ{\sq T} = (\zetaQ{T})_{\Z(\sq
T)}
\end{equation}
to $\Z(\sq T)$ of the central character $\zetaQ{T}$ of $\psi^G$.
\end{subequations}
Hence the factor group $G/\Ker(T)$ is the central product of its subgroup $\sq G\Ker(T)/\Ker(T) \simeq \sq G/\Ker(\sq T)$ with $\Z(T)/\Ker(T) \le \Z(G/\Ker(T))$.  This has a number of consequences, which we shall describe in terms of $G$ instead of $G/\Ker(T)$.

The equations \eqref{CoverG} and \eqref{CoverZ} imply that inclusion $\sq G
\hookrightarrow G$ induces an isomorphism
\begin{subequations} \mylabel{ITT}
\begin{equation} \mylabel{ITTGZs} i^T_{\sq T} \colon \sq G/\Z(\sq T)
\overset{\sim}{\longrightarrow} G/\Z(T)
\end{equation}
of factor groups, sending any coset $\bar{\sigma} \in \sq G/\Z(\sq T)$ to
its product 
\begin{equation} \mylabel{ITTDef} i^T_{\sq T}(\bar{\sigma}) =
\bar{\sigma}\Z(T) \in G/\Z(T)
\end{equation}
\end{subequations}
with $\Z(T)$. In order to exploit this isomorphism we define $\ZG(T)$ to be
the family
\begin{equation} \mylabel{ZG} \ZG(T) = \{\, H \mid \Z(T) \le H \le G\,\}
\end{equation}
of all subgroups of $G$ containing $\Z(T)$. Of course, this family is
defined for any triple $T \in \mfrt$, so that $\ZG(\sq T)$ is the family of
all subgroups of $\sq G$ containing $\Z(\sq T)$. The subgroups of the
factor group $G/\Z(T)$ are just the $H/\Z(T)$ for $H \in \ZG(T)$. Under the
isomorphism $i^T_{\sq T}$ these correspond one to one to the subgroups $\sq
H/\Z(\sq T)$ of $\sq G/\Z(T)$, for $\sq H \in \ZG(\sq T)$. The resulting
bijection of $\ZG(\sq T)$ onto $\ZG(T)$ sends any $\sq H \in \ZG(\sq T)$ to
\begin{subequations} \mylabel{CoverHCorr}
\begin{equation}
H = \sq H\Z(T) \in \ZG(T)  \mylabel{CoverH}
\end{equation}
The inverse bijection sends any $H \in \ZG(T)$ to
\begin{equation}
\sq H = \sq G \cap H \in \ZG(\sq T). \mylabel{CoverHPrime}
\end{equation}
\end{subequations}
Clearly the latter bijection sends $G \in \ZG(T)$ to $\sq G \in \ZG(\sq
T)$, and $N \in \ZG(T)$ to $\sq N = \sq G \cap N \in \ZG(\sq T)$. In view
of \eqref{CoverZ} it also sends $\Z(T) \in \ZG(T)$ to $\Z(\sq T) = \sq G
\cap \Z(T) \in \ZG(\sq T)$.

If $\sq H \in \ZG(\sq T)$ corresponds to $H = \sq H\Z(T) \in \ZG(T)$, then
we may apply Proposition \ref{HKProp} with $\Z(T)$ and its $G$-invariant
linear character $\zetaQ{T}$ in place of $K$ and $\kappa$, respectively.
Since $\zetaQ{\sq T}$ is the restriction \eqref{CoverZeta} of $\zetaQ{T}$
to $\Z(\sq T) = \sq H \cap \Z(T)$, that proposition gives us a bijection of
$\Irr(\,H \mid \zetaQ{T}\,)$ onto $\Irr(\,\sq H \mid \zetaQ{\sq T}\,)$,
sending any $\phi \in \Irr(\,H \mid \zetaQ{T}\,)$ to its restriction
\begin{subequations} \mylabel{CoverPhiCorr}
\begin{equation}
\sq{\phi} = \phi_{\sq H} \in \Irr(\,\sq H \mid \zetaQ{\sq T}\,).
\mylabel{CoverPhiPrime} 
\end{equation}
The inverse bijection sends any $\sq{\phi} \in \Irr(\,\sq H \mid
\zetaQ{\sq T}\,)$ to the character
\begin{equation} \phi = \sq{\phi} * \zetaQ{T} \in \Irr(\,H \mid
\zetaQ{T}\,), \mylabel{CoverPhiZeta}
\end{equation}
whose value at $\sigma\tau \in \sq H\Z(T) = H$ is
\begin{equation} (\sq{\phi} *\zetaQ{T})(\sigma\tau) =
\sq{\phi}(\sigma)\zetaQ{T}(\tau) \in \CC,  \mylabel{CoverPhi}
\end{equation}
\end{subequations}
for any $\sigma \in \sq H$ and $\tau \in \Z(T)$. It follows from
\eqref{CoverT} that the former bijection for $H = N$ sends $\psi \in
\Irr(\,N \mid \zetaQ{T}\,)$ to $\sq{\psi} = \psi_{\sq N} \in \Irr(\,\sq N
\mid \zetaQ{\sq T}\,)$. That same bijection for $H = \Z(T)$  sends
$\zetaQ{T}$ in $\Irr(\,\Z(T) \mid \zetaQ{T}\,)$ to $\zetaQ{\sq T} =
(\zetaQ{T})_{\Z(\sq T)} \in \Irr(\, \Z(\sq T) \mid \zetaQ{\sq T}\,)$ by
\eqref{CoverZeta}.

Characters which correspond in \eqref{CoverPhiCorr} have stabilizers which
correspond in \eqref{CoverHCorr}.
\begin{proposition} \mylabel{CoverStabProp} Suppose that\/ $H \in \ZG(T)$
corresponds to $\sq H \in \ZG(\sq T)$ in \eqref{CoverHCorr}, and
that\/ $\phi \in \Irr(\,H \mid \zetaQ{T}\,)$ corresponds to $\sq{\phi} \in
\Irr(\, \sq H \mid \zetaQ{\sq T}\,)$ in \eqref{CoverPhiCorr}.
Then the stabilizer $G(\phi)$ belongs to $\ZG(T)$, and corresponds to $\sq
G(\sq{\phi}) = \sq G \cap G(\phi) \in \ZG(\sq T)$ in \eqref{CoverHCorr}.
\end{proposition}
\begin{proof}  The stabilizer $G(\phi)$ is a subgroup of $G$ containing the
domain $H$ of $\phi$. Hence it contains $\Z(T) \le H$. Therefore $G(\phi)$
belongs to $\ZG(T)$. 

Any element $\sigma \in \sq G \cap G(\phi)$ normalizes both $\sq G$ and $H
= \Dom(\phi)$, and so normalizes $\sq H = \sq G \cap H$. Because $\sigma$
fixes $\phi \in \Irr(H)$, and normalizes $\sq H \le H$, it fixes the
restriction $\sq{\phi} = \phi_{\sq H}$. Thus any $\sigma \in \sq G \cap
G(\phi)$ lies in $\sq G(\sq{\phi})$.

Conversely, any element $\tau \in \sq G(\sq{\phi})$ belongs to $\sq G$ and
normalizes $\sq H = \Dom(\sq{\phi})$. So it normalizes the product $H = \sq
H\Z(T)$ of $\sq H$ with $\Z(T) \normaleq G$. Because $\tau$ fixes both
$\sq{\phi}$ and the $G$-invariant character $\zetaQ{T}$, it fixes $\phi =
\sq{\phi} * \zetaQ{T}$ by \eqref{CoverPhi}. Thus any $\tau \in \sq
G(\sq{\phi})$ lies in $\sq G \cap G(\phi)$, and the proposition is
proved.
\end{proof}

We may apply the preceding proposition with $H = N \in \ZG(T)$ and $\phi =
\psi \in \Irr(\,N \mid \zetaQ{T}\,)$. Then $\sq H = \sq G \cap N = \sq N$
and $\sq{\phi} = \psi_{\sq N} = \sq{\psi}$. So the above proposition tells
us that $G(\psi) \in \ZG(T)$ corresponds in \eqref{CoverHCorr} to $\sq
G(\sq{\psi}) \in \ZG(\sq T)$. In this case \eqref{CoverH} is
\begin{subequations} \mylabel{CoverGPsiEqs}
\begin{equation} \mylabel{CoverGPsiZT} G(\psi) = \sq G(\sq{\psi})\Z(T).
\end{equation}
Since $\Z(T) \le N \le G(\psi)$, it follows that
\begin{equation} \mylabel{CoverGPsiN}G(\psi) = \sq G(\sq{\psi})N.
\end{equation}
\end{subequations}

We define $\ZST(T)$ to be the family of all subtriples $T' =
(G',N',\psi')$ of $T$ such that
\begin{subequations} \mylabel{ZSTEqs}
\begin{equation} \mylabel{ZSTG} \Z(T) \le G'.
\end{equation}
Because $N$ contains $\Z(T)$, this inclusion implies that
\begin{equation} 
\mylabel{ZSTN} \Z(T) \le  G' \cap N = N'.
\end{equation}
In view of Proposition \ref{SubTZZetaProp} it also implies that
\begin{equation} \mylabel{ZSTZZeta} \Z(T) \le \Z(T') \quad \text{and} \quad
\zetaQ{T} = (\zetaQ{T'})_{\Z(T)}.
\end{equation}
Since $\psi'$ restricts to a multiple $\psi'(1)\zetaQ{T'}$ of the linear
character $\zetaQ{T'}$ on $\Z(T')$, it follows that
\begin{equation} \mylabel{ZSTPsi} (\psi')_{\Z(T)} = \psi'(1)\zetaQ{T}.
\end{equation}
\end{subequations}

Of course the above definition applies with any triple in $\mfrt$ in place
of $T$. In particular, $\ZST(\sq T)$ is the family of all subtriples
$\sq T'$ of $\sq T$ such that the ambient group $\sq G'$ in $\sq T'$
contains $\Z(\sq T)$.
\begin{proposition}\mylabel{CoverZSTProp} There is a one to one
correspondence between all subtriples $T' \in \ZST(T)$ and all subtriples
$\sq T' \in \ZST(\sq T)$. Here two such $T'= (G',N',\psi')$ and
$\sq T' = (\sq G', \sq N',\sq{\psi}')$ correspond if and
only if the subgroups $G', N' \in \ZG(T)$ correspond respectively to $\sq
G', \sq N' \in \ZG(\sq T)$ in \eqref{CoverHCorr}, while the
character $\psi' \in \Irr(\,N' \mid \zetaQ{T}\,)$ corresponds to
$\sq{\psi}' \in \Irr(\,\sq N' \mid \zetaQ{\sq T}\,)$ in
\eqref{CoverPhiCorr}.
\end{proposition}
\begin{proof} If $T' = (G', N',\psi')$ lies in $\ZST(T)$, then \eqref{ZSTG}
implies that $G' \le G$ belongs to the family $\ZG(T)$ in $\eqref{ZG}$. So
$G'$ corresponds to $\sq G \cap G' \in \ZG(\sq T)$ in \eqref{CoverHCorr}.
Similarly, \eqref{ZSTN} implies that $N' \normaleq G'$ belongs to $\ZG(T)$,
and corresponds to $\sq G \cap N' \normaleq \sq G \cap G'$ in
\eqref{CoverHCorr}. The character $\psi'$ lies in $\Irr(\,N' \mid
\zetaQ{T}\,)$ by \eqref{ZSTPsi}, and corresponds to $(\psi')_{\sq G \cap
N'} \in \Irr(\,\sq G \cap N' \mid \zetaQ{\sq T}\,)$ in
\eqref{CoverPhiCorr}.  Hence the triple $\sq T' = (\sq G \cap G', \sq G
\cap N', \psi_{\sq G \cap N'})$ lies in $\mfrt$.  

Evidently $\sq G \cap G'$ is a subgroup of $\sq G$. Because $\sq N$ is $\sq
G \cap N$, and $N'$ is $G' \cap N$, the intersection $(\sq G \cap G') \cap
\sq N$ is equal to $\sq G \cap N'$. Since $\psi' \in \Irr(N')$ lies under
$\psi \in \Irr(N)$, its restriction $(\psi')_{\sq G \cap N'} \in \Irr(\sq G
\cap N')$ lies under $\sq{\psi} = \psi_{\sq N} \in \Irr(\sq N)$. Therefore
$\sq T'$ is a subtriple of $\sq T = (\sq G, \sq N, \sq{\psi})$. This
subtriple belongs to $\ZST(\sq T)$, since its ambient group $\sq G \cap G'
\in \ZG(\sq T)$ contains $\Z(\sq T) = \sq G \cap \Z(T)$ by \eqref{CoverZ} and \eqref{ZSTG}.  Furthermore, it
corresponds to $T'$ in the sense of the present proposition.  Thus any
triple $T' \in \ZST(T)$ corresponds, in the present proposition, to a
unique triple $\sq T' \in \ZST(\sq T)$.

Now let $\sq T' = (\sq G',\sq N', \sq{\psi}')$ be an arbitrary triple in
$\ZST(\sq T)$.   Since $\sq G' \le \sq G$ contains $\Z(\sq T)$, it lies in
$\ZG(\sq T)$, and corresponds to $\sq G'\Z(T) \in \ZG(T)$ in
\eqref{CoverHCorr}. The normal subgroup $\sq N' \normaleq \sq G'$ contains
$\Z(\sq T)$ by \eqref{ZSTN} for $\sq T' \in \ZST(\sq T)$. Hence it lies in
$\ZG(\sq T)$, and corresponds to the normal subgroup $\sq N'\Z(T) \normaleq
\sq G'\Z(T)$ in \eqref{CoverHCorr}.  The character $\sq{\psi}'$ lies in
$\Irr(\,\sq N' \mid \zetaQ{\sq T}\,)$ by \eqref{ZSTPsi} for $\sq T' \in \ZST(\sq T)$. It
corresponds to $\sq{\psi}'*\zetaQ{T} \in \Irr(\,\sq N'\Z(T) \mid \zetaQ{T}\,)$ in
\eqref{CoverPhiCorr}.  So the triple $T' = (\sq G'\Z(T), \sq N'\Z(T), \psi'
* \zetaQ{T})$ lies in $\mfrt$.

The product $\sq G'\Z(T)$ is a subgroup of $G$. Because $\sq N \in \ZG(\sq
T)$ corresponds to $N = \sq N\Z(T)$ in \eqref{CoverHCorr}, the intersection
$\sq N' = \sq G' \cap \sq N \normaleq \sq G'$ has $\sq G'\Z(T) \cap N$ as
its correspondent $\sq N'\Z(T)$. Since $\sq{\psi'}$ lies under $\sq{\psi}$,
it is clear from \eqref{CoverPhi} that $\sq{\psi'} * \zetaQ{T}$ lies under
the correspondent $\psi = \sq{\psi} * \zetaQ{T}$ of $\sq{\psi}$. Thus $T'$
is a subtriple of $T$. This subtriple lies in $\ZST(T)$, since its ambient
group $\sq G'\Z(T)$ contains $\Z(T)$. By construction it corresponds to
$T'$ in the sense of the present proposition. Therefore any triple $\sq T'
\in \ZST(\sq T)$ corresponds, in the present proposition, to a unique
triple $T' \in \ZST(T)$, and the proof of the proposition is complete.
\end{proof}

Notice that the correspondence in the above proposition satisfies
\begin{lemma} \mylabel{CoverZSTLemma} If\/ $T' \in \ZST(T)$ corresponds to
$\sq T' \in \ZST(\sq T)$ in Proposition \ref{CoverZSTProp}, then $\sq T'$
is a subtriple covering $T'$ modulo $\Z(T')$. Furthermore, $\Z(T')$ lies in
$\ZG(T)$, and corresponds to $\Z(\sq T') \in \ZG(\sq T)$ in
\eqref{CoverHCorr}. Finally, $\zetaQ{T'}$ lies in $\Irr(\,\Z(T') \mid
\zetaQ{T}\,)$, and corresponds to $\zetaQ{\sq T'} \in \Irr(\,\Z(\sq T')\mid
\zetaQ{\sq T}\,)$ in \eqref{CoverPhiCorr}.
\end{lemma}
\begin{proof} Let $T'$ be $(G',N',\psi')$, and $\sq T'$ be $(\sq G', \sq
N', \sq{\psi}')$. The correspondent $\sq G' = \sq G \cap G'$ of $G'$ in
\eqref{CoverHCorr} is a subgroup of $G'$.  The correspondent $\sq N' = \sq
G \cap N'$ of $N' = G' \cap N$ is the intersection $\sq G' \cap \sq N$ of
$\sq G' = \sq G \cap G'$ with $\sq N = \sq G \cap N$. The correspondent
$\sq{\psi}' = (\psi')_{\sq N'} \in \Irr(\sq N')$ of $\psi' \in \Irr(N')$ in
\eqref{CoverPhiCorr} lies under $\psi'$. Hence $\sq T'$ is a subtriple of
$T'$. Evidently $\sq G'$ covers its correspondent $G' = \sq G'\Z(T)$ modulo
$\Z(T)$. Because $\Z(T) \le \Z(T')$ by \eqref{ZSTZZeta}, this implies that
$\sq G'$ covers $G'$ modulo $\Z(T')$. Therefore the subtriple $\sq T'$
covers $T'$ modulo $\Z(T')$.

We know from \eqref{ZSTZZeta} that the subgroup $\Z(T')$ of $G' \le G$
contains $\Z(T)$, and hence belongs to the family $\ZG(T)$ in \eqref{ZG}.
The equation \eqref{CoverZ}, for the subtriple $\sq T'$ covering $T'$
modulo $\Z(T')$, implies that
\[ \Z(\sq T') = \sq G' \cap \Z(T') = (\sq G \cap G') \cap \Z(T') = \sq G
\cap \Z(T'). \]
Thus $\Z(T') \in \ZG(T)$ corresponds to $\Z(\sq T') \in \ZG(\sq T)$ in
\eqref{CoverHCorr}.

By \eqref{ZSTZZeta} the linear character $\zetaQ{T'}$ extends $\zetaQ{T}$,
and hence lies in $\Irr(\,\Z(T') \mid \zetaQ{T}\,)$.  Since $\sq T'$ covers
$T'$ modulo $\Z(T')$, we know from \eqref{CoverZeta} that $\zetaQ{\sq T'}$
is the restriction $(\zetaQ{T'})_{\Z(\sq T')}$, which is the
correspondent in $\Irr(\,\Z(\sq T') \mid \zetaQ{\sq T}\,)$  of $\zetaQ{T'}$
in \eqref{CoverPhiCorr}. Therefore the lemma holds.
\end{proof}

Now we can show that the correpondence of subtriples in Proposition
\ref{CoverZSTProp} preserves linear reductions.
\begin{proposition} \mylabel{CoverLRProp} If\/ $T' \in \ZST(T)$ corresponds
to $\sq T' \in \ZST(\sq T)$ in Proposition \ref{CoverZSTProp}, then any
linear reduction of\/ $T'$ lies in $\ZST(T)$ and corresponds to some linear
reduction of\/ $\sq T'$. Similarly, any linear reduction of\/ $\sq T'$ lies
in $\ZST(\sq T)$ and corresponds to some linear reduction of\/ $T'$.
\end{proposition}
\begin{proof} Let $T'$ be $(G', N,\psi')$, and $\sq T'$ be $(\sq G', \sq
N', \sq{\psi}')$. By \eqref{LinRedQuad} and Proposition
\ref{ZZetaLLambdaProp} any linear reduction of $T'$ has the form
\[ T'(\lambda') = (G'(\lambda'), N'(\lambda'), (\psi')_{\lambda'}), \]
for some $L'$ and $\lambda'$ satisfying
\[ L' \normaleq G', \quad \Z(T') \le L' \le N' \quad \text{and} \quad
\zetaQ{T'} \le \lambda' \in \Lin(\,\psi' \mid L'\,).  \]
The stabilizer $G'(\lambda')$ contains the domain $L'$ of $\lambda'$. So it
contains $\Z(T') \le L'$. Since $\Z(T) \le \Z(T')$ by \eqref{ZSTZZeta},
this implies that $\Z(T) \le G'(\lambda')$. Hence the subtriple
$T'(\lambda')$ of both $T'$ and $T$ lies in $\ZST(T)$.

All the subgroups $G'$, $N'$ and $\Z(T')$ of $G$ contain $\Z(T)$ by
\eqref{ZSTEqs}. So they all lie in $\ZG(T)$. Since $T'$ corresponds to $\sq
T'$ in Proposition \ref{CoverZSTProp}, the subgroups $G'$ and $N'$
correspond to $\sq G' = \sq G \cap G' \in \ZST(\sq T)$ and $\sq N' = \sq G
\cap N' \in \ZST(\sq T)$, respectively, in \eqref{CoverHCorr}.
Furthermore, $\Z(T')$ corresponds to $\Z(\sq T') = \sq G \cap \Z(T')$ by
Lemma \ref{CoverZSTLemma}. It follows that the normal subgroup $L'$ of $G'$
satisfying $\Z(T') \le L' \le N'$ lies in $\ZG(T)$, and corresponds in
\eqref{CoverHCorr} to a normal subgroup $\sq L' = \sq G \cap
L' \in \ZG(\sq T)$ of $\sq G'$ satisfying $\Z(\sq T') \le \sq L' \le \sq
N'$.

The linear character $\lambda'$ of $L'$ extends $\zetaQ{T'}$, and hence
extends $\zetaQ{T} = (\zetaQ{T'})_{\Z(T)}$ (see \eqref{ZSTZZeta}). So
$\lambda'$ belongs to $\Irr(\,L' \mid \zetaQ{T} \,)$.  Its correspondent in 
\eqref{CoverPhiCorr} is its restriction to a linear character
$\sq{\lambda}' = (\lambda')_{\sq L'} \in \Irr(\,\sq L' \mid \zetaQ{\sq
T}\,)$ of $\sq L'$. The restriction $\sq{\lambda}'$ lies under the
correspondent $\sq{\psi}' = (\psi')_{\sq N'}$ of $\psi' \in \Irr(N')$,
since $\lambda'$ lies under $\psi'$. Thus $\sq L'$ and $\sq{\lambda}'$
satisfy the equivalent of \eqref{LLambdaConds} for the triple $\sq T'$.
Hence the linear reduction
\[ \sq T'(\sq{\lambda}') = (\sq G'(\sq{\lambda}'), \sq N'(\sq{\lambda}'),
(\sq{\psi}')_{\sq{\lambda}'}) \]
of $\sq T'$ is defined.  

Proposition \ref{CoverStabProp} for $H = L'$ and $\phi = \lambda'$ tells us
that $G(\lambda')$ belongs to $\ZG(T)$ and corresponds to $\sq
G(\sq{\lambda}') \in \ZG(\sq T)$ in \eqref{CoverPhiCorr}.  Since $G'$ and
$N'$ correspond to $\sq G'$ and $\sq N'$, respectively, it follows that
$G'(\lambda') = G' \cap G(\lambda')$ corresponds to $\sq G'(\sq{\lambda'})
= \sq G' \cap \sq G(\sq{\lambda}')$, and that $N'(\lambda') = N' \cap
G(\lambda')$ corresponds to $\sq N'(\sq{\lambda}') = \sq N' \cap \sq
G(\sq{\lambda}')$. Because $(\psi')_{\lambda'} \in \Irr(N'(\lambda'))$ lies
over $\lambda' \in \Lin(\,L' \mid \zetaQ{T}\,)$ and under $\psi' \in
\Irr(\,N' \mid \zetaQ{T}\,)$, it belongs to $\Irr(\,N' \mid \zetaQ{T}\,)$,
and corresponds in \eqref{CoverPhiCorr} to the unique irreducible character
$(\sq{\psi}')_{\sq{\lambda}'}$ of $\sq G'(\sq{\lambda}')$ lying over the
correspondent $\sq{\lambda}' \in \Irr(\,\sq L' \mid \zetaQ{\sq T}\,)$ of
$\lambda'$, and under the correspondent $\sq{\psi}' \in \Irr(\, \sq N' \mid
\zetaQ{\sq T}\,)$ of $\psi'$. We conclude that $T'(\lambda') \in \ZST(T)$
corresponds to $\sq T'(\sq{\lambda}') \in \ZST(\sq T)$ in Proposition
\ref{CoverZSTProp}. Thus the first statement of the present proposition is
proved. The remaining statement is proved similarly, using the fact that
all our correspondences are one to one.
\end{proof}

Our correspondences send multilinear reductions to multilinear
reductions.
\begin{proposition}\mylabel{CoverMLRProp} Any multilinear reduction
of\/ $T$ lies in $\ZST(T)$, and corresponds in Proposition
\ref{CoverZSTProp} to some multilinear reduction of\/ $\sq T$.
Similarly, any multilinear reduction of\/ $\sq T$ lies in $\ZST(\sq
T)$, and corresponds to some multilinear reduction of\/ $T$.
\end{proposition}
\begin{proof} By definition any multilinear reduction $T'$ of $T$ is
the last triple $T' = T_n$ in a finite linear reducing chain $T_0 = T,
T_1, \dots,T_n$ of subtriples of $T$, where each $T_i$, for $i =
1,2,\dots,n$, is a linear reduction of its predecessor $T_{i-1}$. The
initial triple $T_0 = T$ in this chain lies in $\ZST(T)$, and its
correspondent $\sq T_0$ in Proposition \ref{CoverZSTProp} is precisely $\sq
T$. Suppose that $T_{i-1}$ lies in $\ZST(T)$, for some $i = 1,2,\dots,n$, and corresponds to $\sq
T_{i-1} \in \ZST(\sq T)$.  Then Proposition
\ref{CoverLRProp} says that the linear reduction $T_i$ of $T_{i-1}$ lies in
$\ZST(T)$, and corresponds to a linear reduction $\sq T_i \in \ZST(\sq
T)$ of $\sq T_{i-1}$.  By induction we deduce that all the $T_i$, for $i =
0,1,\dots,n$, lie in $\ZST(T)$, and that their correspondents $\sq T_i$ in
$\ZST(\sq T)$ form a linear reducing chain $\sq T_0 = \sq T, \sq T_1,
\dots, \sq T_n$, whose last triple $\sq T_n$ is the correspondent $\sq T'$
of $T' = T_n \in \ZST(T)$. Thus $\sq T'$ is a multilinear reduction of
$T'$, and the first statement of the proposition is proved. The other
statement is proved similarly.
\end{proof}

Finally, our correspondences send linear limits to linear limits. Since
this is the result we shall need in the next section, and its proof has
involved many steps, we state it as
\begin{theorem} \mylabel{CoverLLThm} Suppose $\sq T$ is a subtriple
covering some $T \in \mfrt$ modulo $\Z(T)$. Then the linear limits of\/
$T$ lie in $\ZST(T)$, and correspond one to one in Proposition
\ref{CoverZSTProp} to the linear limits of\/ $\sq T$. If\/ $T' \in
\LL(T)$ corresponds to $\sq T' \in \LL(\sq T)$ in this fashion, then $\sq
T'$ is a subtriple covering $T'$ modulo $\Z(T')$.
\end{theorem}
\begin{proof} Any linear limit $T'$ of $T$ is a multilinear reduction
of $T$. Proposition \ref{CoverMLRProp} tells us that $T'$ lies in
$\ZST(T)$, and that its correspondent $\sq T' \in \ZST(\sq T)$ in
Proposition \ref{CoverZSTProp} is a multilinear reduction of $\sq T$.
By Proposition \ref{CoverLRProp} any linear reduction $\sq T''$ of $\sq T'$
lies in $\ZST(\sq T)$, and corresponds in Proposition \ref{CoverZSTProp} to
some linear reduction $T''$ of $T'$. Any such $T''$ must equal the linear
limit $T'$ of $T$. Hence any linear reduction $\sq T''$ of $\sq T'$ must
equal $\sq T'$. So $\sq T'$ is a linear limit of $\sq T$.

The above argument shows that any linear limit of $T$ lies in $\ZST(T)$,
and corresponds in Proposition \ref{CoverZSTProp} to some linear limit of
$\sq T$. A similar argument shows that any linear limit of $\sq T$ lies in
$\ZST(\sq T)$, and corresponds to some linear limit of $T$. Thus the linear
limits of $T$ correspond one to one to those of $\sq T$ in Proposition
\ref{CoverZSTProp}. The remaining statement in the theorem is a special
case of the first statement in Lemma \ref{CoverZSTLemma}.
\end{proof}

\section{ Equivalence } \mylabel{Equiv}

We say that two subtriples $T'$ and $T''$ of our arbitrary triple $T
= (G,N,\psi) \in \mfrt$ are \emph{equivalent} (and write $T' \sim T''$) if
there is some finite chain $T_0, T_1, \dots,T_n$ of subtriples $T_j \le T$,
starting with $T_0 = T'$ and ending with $T_n = T''$, such that, for each
$j = 1,2,\dots,n$, at least one of the following conditions is satisfied:
\begin{subequations} \mylabel{EquivConds}
\begin{flalign} &\text{$(T_{j-1})^{\tau} = T_j$, for some $\tau \in N$, or}
& & \mylabel{EquivConj} \\
&\text{$T_{j-1}$ is a subtriple covering $T_j$ modulo $\Z(T_j)$, or} & &
\mylabel{EquivCovers} \\
&\text{$T_j$ is a subtriple covering $T_{j-1}$ modulo $\Z(T_{j-1})$.} & &
\mylabel{EquivCovered}
\end{flalign}
\end{subequations}
We call any such chain an \emph{equivalence chain} from $T'$ to $T''$. 
Clearly $\sim$ is an equivalence relation among subtriples of
$T$. 

The notion of equivalent subtriples depends upon $T$, since each member
$T_j$ in the above equivalence chain must be a subtriple of $T$, and the
conjugation in \eqref{EquivConj} is by an element $\tau$ of the normal
subgroup $N$ in $T$. We sometimes speak of \emph{$T$-equivalence} and
\emph{$T$-equivalence chains} to emphasize this dependence.  We do so in
\begin{proposition} \mylabel{STEquivProp} If\/ $\sq T = (\sq G, \sq N,
\sq{\psi})$ is any subtriple of\/ $T$, then any $\sq T$-equivalent
subtriples $T'$ and\/ $T''$ of\/ $\sq T$ are also $T$-equivalent subtriples
of\/ $T$.
\end{proposition}
\begin{proof} There is some $\sq T$-equivalence chain $T_0, T_1, \dots,
T_n$ from $T'$ to $T''$. To prove the proposition we need only show that
these $T_j$ also form a $T$-equivalence chain from $T'$ to $T''$.  Each
$T_j$ is a subtriple of $T$, since $T_j \le \sq T$ and $\sq T \le T$.  The
chain starts with $T_0 = T'$, and ends with $T_n = T''$. For each $j =
1,2,\dots,n$, one of the conditions \eqref{EquivConds} holds with $\sq T$
in place of $T$. Since $\sq N$ is a subgroup of $N$, the condition
\eqref{EquivConj} for $\sq T$ implies the same condition for $T$. The other
two conditions \eqref{EquivCovers} and \eqref{EquivCovered} do not depend
on $T$. So $T_0, T_1, \dots, T_n$ is a $T$-equivalence chain from $T'$ to
$T''$, and the proposition is proved.
\end{proof}

As a general rule, however, we only speak of equivalence and equivalence
chains, with $T$ being understood, as in
\begin{proposition} \mylabel{LLEquivProp} If two subtriples $T'$ and\/
$T''$ of\/ $T$ are equivalent, then any linear limit $\sq T'$ of\/ $T'$ is
equivalent to some linear limit $\sq T''$ of\/ $T''$.
\end{proposition}
\begin{proof} We use induction on the length $n \ge 0$ of an equivalence
chain $T_0, T_1, \dots, T_n$ from $T'$ to $T''$. If $n = 0$, then $T' = T_0
= T_n = T''$, and the proposition is trivial. So we may assume that $n >
0$, and that the proposition holds with $T_{n-1}$ in place of $T'' = T_n$.
Then any linear limit $\sq T'$ of $T'$ is equivalent to some linear limit
$\sq T_{n-1}$ of $T_{n-1}$.

One of the three conditions \eqref{EquivConds} holds with $j = n$. If
\eqref{EquivConj} holds, then $(T_{n-1})^{\tau} = T_n$, for some $\tau \in
N$. Proposition \ref{ConjLLProp} tells us that conjugation by $\tau$ sends
$\sq T_{n-1} \in \LL(T_{n-1})$ to some linear limit $\sq T'' = (\sq
T_{n-1})^{\tau}$ of $T'' = T_n$. If \eqref{EquivCovers} holds, then the
subtriple $T_{n-1}$ covers $T_n$ modulo $\Z(T_n)$. In that case Theorem
\ref{CoverLLThm} says that $\sq T_{n-1} \in \LL(T_{n-1})$ is a subtriple
covering some linear limit $\sq T''$ of $T'' = T_n$ modulo $\Z(\sq T'')$.
That theorem also applies when \eqref{EquivCovered} holds, i.~e., when the
subtriple $T_n$ covers $T_{n-1}$ modulo $\Z(T_{n-1})$. In that case it says
that some linear limit $\sq T''$ of $T'' =T_n$ is a subtriple covering $\sq
T_{n-1} \in \LL(T_{n-1})$ modulo $\Z(\sq T_{n-1})$. In each of the three
cases $\sq T_{n-1}$ is equivalent to a linear limit $\sq T''$ of $T''$. So
$\sq T' \sim \sq T_{n-1}$ is equivalent to $\sq T''$, and the proposition
is proved.
\end{proof}

In the language of the present section the Main Theorem in the introduction
can be stated as
\begin{theorem} \mylabel{LLThm} All linear limits of any triple $T =
(G,N,\psi) \in \mfrt$ are equivalent.
\end{theorem}
\begin{proof}  We prove this theorem by induction on the order $|G|$ of the
finite group $G$.  If $|G| = 1$, then the only possible linear
limit of $T$ is $T$ itself. In this case the theorem is trivial. So we may
assume that $|G| > 1$, and that the theorem holds for all strictly smaller
values of $|G|$. We may also assume that the theorem does not hold for $T$.
We divide the rest of the proof into a number of steps, all based on these
assumptions. These steps will lead to a contradiction, thus proving the
theorem.

\begin{step} \mylabel{TNotLimit} $T \notin \LL(T)$.
\end{step}
\begin{proof} If $T \in \LL(T)$, then the only possible linear reduction of
$T$ is $T$ itself. So $T$ is the only member of $\LL(T)$. Thus the theorem
holds trivially in this case, contradicting our assumptions. This
contradiction proves the present step.
\end{proof}

Let $\LSG$ be the family of all subgroups $L$ satisfying
\begin{subequations} \mylabel{LLSubsets}
\begin{equation} \mylabel{GoodL} 
L \normaleq G, \quad  \Z(T) < L \le N,\quad \text{and} \quad \Lin(\,\psi
\mid L \,) \ne \emptyset, 
\end{equation}
where $\emptyset$ is the usual empty set. For each $L \in \LSG$ we
define $\LL(\,T \mid L\,)$ to be the non-empty subset
\begin{equation} \mylabel{LLTL}  \LL(\,T \mid L\,) = \bigcup_{\lambda \in
\Lin(\psi \mid L)} \LL(T(\lambda))
\end{equation}
\end{subequations}
of $\LL(T)$. Then we have
\begin{step} \mylabel{LLUnion}  $\LL(T) = \bigcup_{L \in \LSG} \LL(\,T
\mid L\,)$.  \end{step}
\begin{proof} We only need show that any $T' \in \LL(T)$ belongs to the
subset $\LL(\,T \mid L\,) \subseteq \LL(T)$, for some $L \in
\LSG$. The linear limit $T'$ of $T$ is different from $T$ by Step
\ref{TNotLimit}. It follows that there is some finite linear reducing chain
$T_0, T_1, \dots,T_n$, from $T = T_0$ to $T' = T_n$, such that $n > 0$, and
each $T_i$, for $i = 1,2,\dots,n$, is a proper linear reduction of
$T_{i-1}$. In particular, $T'$ is a linear limit of the proper linear
reduction $T_1$ of $T$. By Propositions \ref{PropLinRedProp} and
\ref{ZZetaLLambdaProp} there exist some $L$ satisfying \eqref{GoodL}, and
some $\lambda \in \Lin(\,\psi \mid L\,)$, such that $T_1 = T(\lambda)$.
Then $T' \in \LL(\,T(\lambda))$ belongs to the subset $\LL(\,T \mid L\,)$
in \eqref{LLTL}, and the present step is proved.
\end{proof}

Our induction hypothesis will give
\begin{step}  \mylabel{LTLEquiv} If $L \in \LSG$, then any $T', T'' \in
\LL(\, T \mid L\,)$ are equivalent.
\end{step}
\begin{proof} By \eqref{LLTL} any triples $T'$ and $T''$ in $ \LL(\,T \mid
L\,)$ are linear limits of $T(\lambda')$ and $T(\lambda'')$, respectively,
for some $\lambda', \lambda'' \in \Lin(\,\psi \mid L\,)$. Since $\psi$ is
an irreducible character of $N$, Clifford theory for $L \normaleq N$ tells
us that $\lambda'' = (\lambda')^{\tau}$, for some $\tau \in N$.  So
$(T')^{\tau}$ is a linear limit of $T(\lambda'') = T((\lambda')^{\tau}) =
T(\lambda')^{\tau}$. By Proposition \ref{PropLinRedProp} the inclusion
$\Z(T) < L$ in \eqref{GoodL} implies that $T(\lambda'')$ is a proper linear
reduction of $T$. Hence the ambient group $G(\lambda'')$ in
$T(\lambda'')$ has order $|G(\lambda'')| < |G|$, by that same proposition.
Our induction hypothesis tells us that the two linear limits $(T')^{\tau}$
and $T''$ of $T(\lambda'')$ are $T(\lambda'')$-equivalent. Proposition
\ref{STEquivProp} implies that they are $T$-equivalent. Therefore $T' \sim
(T')^{\tau} \sim T''$, and this step is proved.
\end{proof}

Another useful observation is
\begin{step} \mylabel{LKLIncl} If\/ $L,K \in \LSG$ satisfy $K \le
L$, then $\LL(\,T \mid L\,) \subseteq \LL(\,T \mid K\,)$.
\end{step}
\begin{proof} By \eqref{LLTL} any $T' \in \LL(\,T \mid L\,)$ is a linear
limit of $T(\lambda)$, for some $\lambda \in \Lin(\,\psi \mid L\,)$.  Then
$\kappa = \lambda_K$ lies in $\Lin(\, \psi \mid K\,)$, so that $T(\kappa)$
is defined. Proposition \ref{CompoundLinRedProp} tells us that $T(\lambda)$
is also the linear reduction $[T(\kappa)](\lambda)$ of $T(\kappa)$. Hence
its linear limit $T'$ is also a linear limit of $T(\kappa)$, and so lies in
$\LL(\,T \mid K\,)$. Thus the present step holds.
\end{proof}

As a consequence of the preceding two steps we have
\begin{step} \mylabel{LMInt} If\/ $\Z(T) < L \cap M$, for some $L,M \in
\LSG$, then $T'\sim T''$, for any $T' \in \LL(\,T \mid L\,)$ and\/ $T''
\in \LL(\,T \mid M\,)$.
\end{step}
\begin{proof} The intersection $K = L \cap M$ is a normal subgroup of $G$
contained in $N$, since both $L$ and $M$ are such subgroups. Any character in the non-empty set $\Lin(\,\psi \mid L\,)$ restricts to
a character in $\Lin(\,\psi \mid K\,)$. Since $\Z(T) < K$ by
assumption, we conclude that $K$ belongs to $\LSG$. Now Step \ref{LKLIncl},
for both $K \le L$ and $K \le M$, implies that both $T' \in \LL(\,T \mid
L\,)$ and $T'' \in \LL(\,T \mid M\,)$ belong to $\LL(\,T \mid K\,)$.  So
$T' \sim T''$ by Step \ref{LTLEquiv}, and the present step is proved.
\end{proof}

Since the theorem is false for $T$, there exist two linear limits $T'$ and
$T''$ of $T$ such that
\begin{subequations} \mylabel{TTConds}
\begin{equation} \mylabel{TTNotEquiv} T' \not \sim T''.
\end{equation}
In view of Step \ref{LLUnion} there are two subgroups $L, M$ satisfying
\begin{equation} \mylabel{LMConds} L, M \in \LSG, \quad T' \in \LL(\,T \mid
L\,), \quad \text{\and}\quad T'' \in \LL(\,T \mid M\,).
\end{equation}
Then \eqref{LLTL} gives us two linear characters $\lambda, \mu$ such that
\begin{equation} \mylabel{LambdaMuConds}
\lambda \in \Lin(\,\psi \mid L\,), \quad \mu \in \Lin(\,\psi \mid M\,),
\quad T' \in \LL(T(\lambda)), \quad \text{and} \quad T'' \in \LL(T(\mu)).
\end{equation}
\end{subequations}
For the rest of the proof of the theorem we fix $T'$, $T''$, $L$, $M$,
$\lambda$, and $\mu$ with these properties.
\begin{step} \mylabel{LambdaMuZeta} The intersection $L \cap M$ is $\Z(T)$.
Both $\lambda \in \Lin(L)$ and $\mu \in \Lin(M)$ restrict to $\zetaQ{T}$ on $\Z(T)$.
\end{step}
\begin{proof} Suppose the $L \cap M \ne \Z(T)$. The subgroup $L \in \LSG$
contains $\Z(T)$ by \eqref{GoodL}. So do $M$ and the intersection $L \cap
M$. Hence $\Z(T) < L \cap M$. By Step \ref{LMInt} this implies that $T' \in
\LL(\,T \mid L\,)$ is equivalent to $T'' \in \LL(\,T \mid M\,)$,
contradicting \eqref{TTNotEquiv}. This contradiction proves the first
statement in this step. The other statement follows immediately from the
last statement in Proposition \ref{ZZetaLLambdaProp}.
\end{proof}

By \eqref{GoodL} both $L$ and $M$ are normal subgroups of $G$ contained in
$N$. Hence their commutator $[L,M]$ is contained in their intersection $L
\cap M$, which we have just seen to be $\Z(T)$. So $L$ and $M$ satisfy
\eqref{CLMConds}, and \eqref{C} defines a $G$-invariant, bilinear function
$c$ from $L \times M$ to $C^{\times}$.
\begin{step} \mylabel{MLPerp} The perpendicular subgroup
\[ L^{\perp} = \{\,\sigma \in M \mid c(L,\sigma) = 1\,\}  \]
to $L$ with respect to the form $c$ is a normal subgroup of\/ $G$ contained
in $M \le N$. It is also the stabilizer $M(\lambda)$ of\/ $\lambda \in
\Lin(L)$ in $M$.

The similar perpendicular subgroup $M^{\perp}$ to $M$ with respect to $c$
is a normal subgroup of\/ $G$ contained in $L \le N$, and equals the
stabilizer $L(\mu)$ of\/ $\mu$ in $L$.
\end{step}
\begin{proof} The bilinearity \eqref{BilinearC} of $c$ implies that
$L^{\perp}$ is a subgroup of $M \le N$. This subgroup is normal in
$G$, since $L$ and $M$ are normal subgroups of $G$, and $c$ is
$G$-invariant by \eqref{GInvC}.

The commutator $[\rho,\sigma] = \rho^{-1}\sigma^{-1}\rho\sigma$ is equal to
$\rho^{-1}\rho^{\sigma}$, for any $\rho \in L$ and $\sigma \in M$. Both
$\rho$ and $\rho^{\sigma}$ lie in $L \normaleq G$. Furthermore, $\lambda
\in \Lin(L)$ extends $\zetaQ{T}$ by Step \ref{LambdaMuZeta}. It follows
that
\[ c(\rho, \sigma) = \zetaQ{T}([\rho, \sigma]) =
\lambda(\rho^{-1}\rho^{\sigma}) =
\lambda(\rho)^{-1}\lambda^{\sigma^{-1}}(\rho). \]
Hence $c(\rho, \sigma) = 1$, for all $\rho \in L$, if and only if
$\lambda^{\sigma^{-1}} = \lambda$. This just says that $L^{\perp}$ is equal
to $M(\lambda)$.  Thus the first paragraph of this step is proved. The
other paragraph is proved similarly.
\end{proof}

Now we can construct two new members of $\LSG$.
\begin{step} \mylabel{MLambdaL} The two products $LM(\lambda)$ and
$L(\mu)M$ belong to $\LSG$.
\end{step}
\begin{proof} The subgroup $L \in \LSG$ is normal in $G$ and contained in
$N$ by \eqref{GoodL}. The subgroup $M(\lambda) = L^{\perp}$ is normal in
$G$ and contained in $N$ by Step \ref{MLPerp}. Hence their product
$LM(\lambda) = M(\lambda)L$ is a normal subgroup of $G$ contained in $N$.
The strict inclusion $\Z(T) < L$ in \eqref{GoodL} implies that $\Z(T) <
LM(\lambda)$. 

We must show that $\Lin(\,\psi \mid LM(\lambda)\,)$ is not empty. In view
of Proposition \ref{IsotropicL} this happens if and only if
$[LM(\lambda),LM(\lambda)]$ is contained in $\Ker(T)$. Because $L$ and
$M(\lambda)$ are normal subgroups of $G$, we have
\[  [LM(\lambda),LM(\lambda)] =
[L,L][L,M(\lambda)][M(\lambda),M(\lambda)]. \]
The factor $[L,L]$ is contained in $\Ker(T)$ by Proposition
\ref{IsotropicL}, since $\Lin(\,\psi \mid L\,)$ has an element $\lambda$.
The linear character $\zetaQ{T}$ sends the factor $[L,M(\lambda)]$ to
$c(L,M(\lambda))$, which is $1$ since $M(\lambda) = L^{\perp}$ by Step
\ref{MLPerp}. Hence $[L,M(\lambda)]$ is contained in $\Ker(T) =
\Ker(\zetaQ{T})$. Finally, the factor $[M(\lambda),M(\lambda)]$ is
contained
in $[M,M]$, which is a subgroup of $\Ker(T)$ by Proposition
\ref{IsotropicL}, since $\mu$ lies in $\Lin(\,\psi \mid M\,)$. Thus each
factor in the above product is contained in $\Ker(T)$, and the proof that
$\Lin(\,\psi \mid LM(\lambda)\,) \ne \emptyset$ is complete.

The above arguments show that $LM(\lambda)$ satisfies all the conditions for
$L$ in \eqref{GoodL}, and hence belongs to $\LSG$. The proof that $L(\mu)M
\in \LSG$ is similar.
\end{proof}

At this point things simplify drastically.
\begin{step} \mylabel{MLambda} Both $M(\lambda)$ and $L(\mu)$ are equal to
$\Z(T)$.
\end{step}
\begin{proof} Clearly $\Z(T) = L \cap M$ is contained
in $M$, and fixes $\lambda \in \Lin(L)$. Hence $\Z(T) \le M(\lambda)$. If $\Z(T) <
M(\lambda)$, then $\Z(T) < M \cap LM(\lambda)$. Because
$LM(\lambda) \in \LSG$ by Step \ref{MLambdaL}, there is some $T''' \in
\LL(\,T \mid M(\lambda)L\,)$.  Step \ref{LMInt} now tells us that $T'' \in
\LL(\,T \mid M\,)$ is equivalent to $T'''$. That step also tells us that
$T'''$ is equivalent to $T' \in \LL(\,T \mid L\,)$, since $\Z(T) < L =
M(\lambda)L \cap L$. Hence $T'' \sim T''' \sim T'$, contradicting
\eqref{TTNotEquiv}. This contradiction proves that $M(\lambda) = \Z(T)$.
The proof that $L(\mu) = \Z(T)$ is similar.
\end{proof}

Now we can construct some covering subtriples.
\begin{step} \mylabel{GLambda} Both $G(\mu)L$ and $G(\lambda)M$ are equal
\begin{subequations}
\begin{equation} \mylabel{GMuL} G(\mu)L = G = G(\lambda)M
\end{equation}
to $G$. Hence $G(\lambda, \mu)L$ is $G(\lambda)$, and
\begin{equation} \mylabel{sqTPrime}  \sq T' = (G(\lambda,\mu),
N(\lambda,\mu), (\psi_{\lambda})_{N(\lambda,\mu)}) 
\end{equation}
is a subtriple covering $T(\lambda)$ modulo $\Z(T(\lambda))$. Similarly,
$G(\mu)$ is $G(\lambda,\mu)M$, and
\begin{equation} \mylabel{sqTDoublePrime} \sq T'' = (G(\lambda,\mu),
N(\lambda,\mu), (\psi_{\mu})_{N(\lambda,\mu)}) 
\end{equation}
\end{subequations}
is a subtriple covering $T(\mu)$ modulo $\Z(T(\mu))$.
\end{step}
\begin{proof} Steps \ref{MLPerp} and \ref{MLambda} imply that $\Z(T) =
L(\mu) = M^{\perp}$, and $\Z(T) = M(\lambda) = L^{\perp}$.
Since $c$ is bilinear, it follows that the factor groups $L/\Z(T)$ and
$M/\Z(T)$ are abelian, and that $c$ induces a non-singular bilinear
pairing $\bar c$ of $(L/\Z(T)) \times (M/\Z(T))$ into
$\CC^{\times}$. So $L/\Z(T)$ and $M/\Z(T)$ are dual finite abelian groups.
In particular, they have the same order
\[  [L : \Z(T)] = [M : \Z(T)]. \]

The linear character $\mu$ of $M$ extends $\zetaQ{T} \in \Lin(\Z(T))$ by
Step \ref{LambdaMuZeta}. Since $M/\Z(T)$ is abelian, it follows that
$\Lin(\,M \mid \zetaQ{T}\,) = \Irr(\,M \mid \zetaQ{T}\,)$ consists of 
$[M : \Z(T)]$ distinct extensions of $\zetaQ{T}$ to linear characters of
$M$. Because $G$ normalizes $M$, and leaves $\zetaQ{T}$ invariant, it
permutes these extensions among themselves by conjugation. So does the
subgroup $L \le G$. The $L$-orbit  of $\mu \in \Lin(\,M \mid \zetaQ{T}\,)$
under this action has order
\[ [L : L(\mu)] = [L : \Z(T)] = [M : \Z(T) ] = \bigl|\Lin(\,M \mid
\zetaQ{T}\,)\bigr |, \]
since $L(\mu) = \Z(T)$ by Step \ref{MLambda}. Hence this $L$-orbit is all
of $\Lin(\,M \, \mid \zetaQ{T}\,)$. In particular, this $L$-orbit is also
the $G$-orbit of $\mu$. This implies  the first equation in \eqref{GMuL}.
The other equation there is proved similarly. Thus the first statement of
this step holds.

The subgroup $L \le G$ fixes its own character $\lambda$, and so is
contained in $G(\lambda)$. This and the first equation in \eqref{GMuL}
imply that $G(\lambda)$ is the product $G(\lambda, \mu)L$ of
$G(\lambda,\mu) = G(\lambda) \cap G(\mu)$ with $L$. The $\lambda$-Clifford
correspondent $\psi_{\lambda} \in\Irr(\,N(\lambda) \mid \lambda\,)$
restricts to a multiple of $\lambda \in \Lin(L)$. So Proposition
\ref{UnderZetaQProp}, with $T(\lambda) = (G(\lambda), N(\lambda),
\psi_{\lambda})$ in place of $T$, tells us that $L \le \Z(T(\lambda))$.
Hence
\[ G(\lambda) = G(\lambda,\mu)L = G(\lambda,\mu)\Z(T(\lambda)). \]

Now we may apply Proposition \ref{CoverTProp}, with $T(\lambda)$ and
$G(\lambda,\mu)$ in place of $T$ and $\sq G$, respectively, to obtain the
second statement of the present step. The final statement is proved
similarly.
\end{proof}

\begin{step} \mylabel{EqualChars} The two irreducible characters
$(\psi_{\lambda})_{N(\lambda,\mu)}$ and\/ $(\psi_{\mu})_{N(\lambda,\mu)}$
of\/ $N(\lambda,\mu)$ are equal.  Hence the two subtriples $\sq T'$ and
$\sq T''$ in Step \ref{GLambda} are equal, and the two linear reductions
$T(\lambda)$ and $T(\mu)$ are equivalent.
\end{step}
\begin{proof} Since $N(\mu)$ contains $M = \Dom(\mu) \le N$, it follows
from the second equation in \eqref{GMuL} that
\[ N = N \cap G(\lambda)M = (N \cap G(\lambda))M = N(\lambda)N(\mu). \]
So $N(\lambda)N(\mu)$ is the only $N(\lambda), N(\mu)$-double coset in $N$,
and the intersection $N(\lambda) \cap N(\mu)$ is $N(\lambda, \mu)$. Because
both $\psi_{\lambda} \in \Irr(N(\lambda)$) and $\psi_{\mu} \in
\Irr(N(\mu))$ induce $\psi \in \Irr(N)$, Mackey's Formula gives the inner
products
\[  1 = \gen{\psi,\psi} = \gen{(\psi_{\lambda})^N, (\psi_{\mu})^N} =
\gen{(\psi_{\lambda})_{N(\lambda,\mu)}, (\psi_{\mu})_{N(\lambda,\mu)}}. \]
But the characters $(\psi_{\lambda})_{N(\lambda,\mu)}$ and
$(\psi_{\mu})_{\N(\lambda,\mu)}$ in the two triples \eqref{sqTPrime} and
\eqref{sqTDoublePrime} are both irreducible characters of
$N(\lambda,\mu)$. So the fact that their inner product is $1$ forces them
to be equal. 

The two triples $\sq T'$ and $\sq T''$ already have the same ambient groups
and normal subgroups. We have just shown that they have the same
characters. Hence they are equal. In view of Step \ref{GLambda},
this implies that the same subtriple $\sq T' = \sq T''$ covers both
$T(\lambda)$ and $T(\mu)$ modulo their centers. So $T(\lambda)$ is
equivalent to $T(\mu)$, and the present step is proved.
\end{proof}

We can finally finish the  proof of Theorem \ref{LLThm}. By Proposition
\ref{LLEquivProp} the equivalence $T(\lambda) \sim T(\mu)$ in Step
\ref{EqualChars} implies that the linear limit $T''$ of $T(\mu)$ in
\eqref{LambdaMuConds} is equivalent to some linear limit $T'''$ of
$T(\lambda)$.  This $T'''\in \LL(T(\lambda))$ belongs to the set $\LL(\,T
\mid L\,)$ in \eqref{LLTL}. So does $T'$ by \eqref{LMConds}. Therefore
$T''' \sim T'$ by Step \ref{LTLEquiv}. Thus we have equivalences $T'' \sim
T''' \sim T'$, which contradict \eqref{TTNotEquiv}. This contradiction
completes the proof of Theorem \ref{LLThm}.
\end{proof}

\section{ Symplectic limits } \mylabel{Symp} 

We're going to examine consequences of various restrictions on the normal
section
\begin{equation} \mylabel{ST} \SSS(T) = N/\Z(T)
\end{equation}
of our arbitrary triple $T = (G,N,\psi) \in \mfrt$, especially when $T$ is
linearly irreducible. Many of these consequences are standard steps in the
Hall-Higman reduction. This is certainly true for
\begin{proposition} \mylabel{NilpSTProp} If\/ $\SSS(T)$ is nilpotent, and\/
$T$ is linearly irreducible, then $\SSS(T)$ is abelian.
\end{proposition}
\begin{proof} Let $M$ be the inverse image in $N$ of the center
$\Z(\SSS(T))$ of the nilpotent group $\SSS(T) = N/\Z(T)$. Then $M$ is a
normal subgroup of $G$ such that
\[  \Z(T) \le M \le N \quad \text{and} \quad M/\Z(T) = \Z(N/\Z(T)). \]
Because $M/\Z(T)$ is central in $N/\Z(T)$, the commutator $[N,M]$ is
contained in $\Z(T)$.  So the conditions \eqref{CLMConds} are satisfied
with $L = N$. Hence  \eqref{C} defines a $G$-invariant, bilinear function
$c$ from $N \times M$ to $\CC^{\times}$. The perpendicular subgroup
$ M^{\perp} = \{\, \rho \in N \mid c(\rho, M) = 1\, \} $
to $M$ with respect to this function is a normal subgroup of $G$ contained
in $N$. Since $c$ is bilinear, the factor group $N/M^{\perp}$ is abelian.
So we can prove the proposition by showing that
\[   M^{\perp} = \Z(T). \]

Suppose this is false. Proposition \ref{CZPerpProp} tells us that $\Z(T) =
\Z(T) \cap N \le M^{\perp}$. Hence $\Z(T) < M^{\perp}$. So
$M^{\perp}/\Z(T)$ is a non-trivial normal subgroup of the nilpotent group
$N/\Z(T)$. It follows that $M^{\perp}/\Z(T)$ has a non-trivial intersection
with $M/\Z(T) = \Z(N/\Z(T))$. We conclude that $L = M^{\perp} \cap M$ is a
normal subgroup of $G$ such that $\Z(T) < L \le N$. Furthermore, $[L,L] \le
[N,M] \le \Z(T)$, and $\zetaQ{T}([L,L]) = c(L,L) \le c(M^{\perp}, M) = 1$.
Thus $[L,L]$ is contained in $\Ker(\zetaQ{T}) = \Ker(T)$ (see
\eqref{KerZetaQ}).  Since $\Z(T) < L$, the last statement in Proposition
\ref{PropLinRedProp} tells us that $T$ is linearly reducible, contrary to
our hypotheses. This contradiction proves the proposition.
\end{proof}

From now on we assume that $\SSS(T)$ is abelian. The group $G$ acts on its
normal section $\SSS(T) = N/\Z(T)$ by conjugation. Since $\SSS(T)$ is
abelian, this induces an action of the factor group $G/N$ as automorphisms
of $\SSS(T)$, with the coset $\tau N \in G/N$, for any $\tau \in G$,
sending any coset $\bar{\sigma} \in \SSS(T) = N/\Z(T)$ to
\begin{equation} \mylabel{GNOnST} \bar{\sigma}^{\tau N} =
\bar{\sigma}^{\tau} \in \SSS(T).
\end{equation}
The commutator $[N,N]$ is contained in $\Z(T)$, because $[N/\Z(T),N/\Z(T)]
= 1$. Hence the conditions \eqref{CLMConds} hold with $L = M = N$. So
\eqref{C} defines a $G$-invariant, bilinear function $c = c_T$ from $N \times N$
to $\CC^{\times}$, sending any $\rho, \sigma \in N$ to
\begin{subequations}\mylabel{SCBarEqs}
\begin{equation} \mylabel{SC} c(\rho, \sigma) = \zetaQ{T}([\rho, \sigma]) \in \CC^{\times}.
\end{equation}
As in \eqref{CBar}, this induces a $G/N$-invariant,
bilinear function $\bar c = \bar c_T$ from $\SSS(T) \times \SSS(T)$ to
$\CC^{\times}$, sending the cosets $\rho \Z(T), \sigma \Z(T) \in \SSS(T) =
N/\Z(T)$ to
\begin{equation} \mylabel{SCBar}  \bar c(\rho \Z(T), \sigma \Z(T)) =
\zetaQ{T}([\rho, \sigma]) \in \CC^{\times},
\end{equation}
for any $\rho, \sigma \in N$.  This last bilinear function is clearly
\emph{strongly alternating}, in the sense that
\begin{equation} \mylabel{StrAltCBar} \bar c(\bar{\rho}, \bar{\rho}) = 1,
\end{equation}
for any $\bar{\rho} \in \SSS(T)$. It follows that it is \emph{alternating},
in the sense that
\begin{equation} \mylabel{AltCBar} \bar c(\bar{\sigma}, \bar{\rho}) = \bar
c(\bar{\rho}, \bar{\sigma})^{-1}, 
\end{equation}
for all $\bar{\rho}, \bar{\sigma} \in \SSS(T)$. 
\end{subequations}

Let $H$ be any finite group. We need a short name for the situation
consisting of a finite abelian $H$-group $A$, together with an
$H$-invariant, strongly alternating, bilinear function $b \colon A \times A
\to \CC^{\times}$.  Since $b$ behaves much like an alternating bilinear
form on a vector space, we shall call such an $A$ and $b$ a \emph{formed
abelian $H$-group}. Usually we just speak of ``the formed abelian $H$-group
$A$,'' with the \emph{bilinear form} $b = b_A$ for $A$ being understood. 
Notice that $\SSS(T)$ is a formed abelian $G/N$-group, with the action
\eqref{GNOnST} and the bilinear form $\bar c$ in \eqref{SCBar}, whenever
$\SSS(T)$ is abelian.

The terminology associated with an arbitrary formed abelian $H$-group $A$
is adapted from that for bilinear forms on a vector space. Because the
bilinear form $b$ for $A$ is alternating, the left and right
\emph{perpendicular subgroups} $B^{\perp}$ to any subgroup $B \le A$ coincide
\begin{subequations} \mylabel{HAEqs}
\begin{equation} \mylabel{PerpGrp} B^{\perp} = \{\,\rho \in A \mid b(\rho,
B) = 1 \, \} = \{\, \sigma \in A \mid b( B,\sigma) = 1\,\}.
\end{equation} 
The $H$-invariance of $b$ implies that $B^{\perp}$ is an $H$-invariant
subgroup of $A$ whenever $B$ is one. In particular, the \emph{radical} 
\begin{equation} \mylabel{RadA} \rad(A) = A^{\perp}
\end{equation}
\end{subequations}
of $A$ is an $H$-invariant subgroup of $A$. 

The form $b$ is \emph{non-singular} if $\rad(A) = 1$. In this case $A$ and
$b$ behave like a vector space with a symplectic form. So we say that $A$
is \emph{symplectic} if $b$ is non-singular, and define a \emph{symplectic
$H$-group} to be a symplectic formed abelian $H$-group. When $A$ is
symplectic, its bilinear form $b \colon A \times A \to \CC^{\times}$ is a
duality of the finite abelian group $A$ with itself. Hence $b$ induces a
duality between any subgroup $B$ of $A$ and the factor group $A/B^{\perp}$.
In particular, these two finite groups have the same order
\begin{subequations} \mylabel{SympEqs}
\begin{equation} \mylabel{BDualOrder} |B| = [A : B^{\perp}].
\end{equation}
It follows that $B^{\perp\perp} = (B^{\perp})^{\perp}$ has the same index
$[A : B^{\perp\perp}] = [A : B]$ in $A$ as its subgroup $B$. Therefore
\begin{equation} \mylabel{BPerpPerp} B^{\perp\perp} = B
\end{equation}
\end{subequations}
when $A$ is symplectic.

A subgroup $B$ of an arbitrary formed abelian $H$-group $A$ is
\emph{isotropic} if $b(B,B) = 1$, i.~e., if $B \le B^{\perp}$. We say that
$A$ is \emph{$H$-anisotropic} if $1$ is the only $H$-invariant isotropic
subgroup of $A$. Since $\rad(A)$ is always an $H$-invariant isotropic
subgroup, any $H$-anisotropic $A$ is a symplectic $H$-group. Furthermore,
any $H$-invariant subgroup $B$ of $A$ is also an $H$-anisotropic formed
abelian $H$-group, with the restriction of $b$ as its bilinear form. The
intersection $B \cap B^{\perp}$ in $A$ is the radical $\rad(B) = 1$ of $B$.
Since $b$ is non-singular, this and \eqref{BDualOrder} imply that $A$ is
the perpendicular direct product
\begin{subequations} \mylabel{AnisEqs}
\begin{equation} \mylabel{AnisComps} A = B \times B^{\perp}
\end{equation}
of its two $H$-anisotropic subgroups $B$ and $B^{\perp}$. By induction it
follows that $A$ is a direct product
\begin{equation} \mylabel{AnisDecomp} A = B_1 \times B_2 \times \dots
\times B_k
\end{equation} 
\end{subequations}
of pairwise orthogonal simple $H$-subgroups $B_1, B_2, \dots, B_k$, for
some integer $k \ge 0$. Notice that each simple $H$-subgroup $B_i$ is an
elementary $p_i$-group, for some prime $p_i$. So the Sylow $p$-subgroups of
$A$ are all elementary abelian $p$-groups when $A$ is anisotropic.

The importance of anisotropy for our formed abelian $G/N$-group $\SSS(T)$
is explained by the next step in the Hall-Higman analysis.
\begin{proposition} \mylabel{AbIrrProp} If\/ $\SSS(T)$ is abelian, then the
$G/N$-invariant isotropic subgroups of\/ $\SSS(T) = N/\Z(T)$ are precisely
the images $L\Z(T)/\Z(T)$ of the normal subgroups $L$ of $G$ contained in
$N$ such that $[L,L] \le \Ker(T)$. Hence $T$ is linearly irreducible if and
only if\/ $\SSS(T)$ is $G/N$-anisotropic.
\end{proposition}
\begin{proof} The $G/N$-invariant subgroups of $\SSS(T) = N/\Z(T)$ are just
the images $\bar L = L\Z(T)/\Z(T)$ of the normal subgroups $L$ of $G$
contained in $N$.  It follows from \eqref{SCBar} that $\bar c(\bar L, \bar
L) = \zetaQ{T}([L,L])$, for any such $\bar L$ and $L$. Hence $\bar L$ is
isotropic if and only if $[L,L]$ is contained in $\Ker(\zetaQ{T})$. This
last group is $\Ker(T)$ by \eqref{KerZetaQ}.  Thus the first statement of
the proposition holds. It implies the remaining statement by Proposition
\ref{PropLinRedProp}.
\end{proof}

When $T$ is linearly irreducible, its $G/N$-anisotropic formed abelian
$G/N$-group $\SSS(T)$ is symplectic. So the following Hall-Higman result applies in
that case.
\begin{proposition} \mylabel{InvPsiProp} If\/ $\SSS(T)$ is a symplectic
$G/N$-group, then $\psi$ vanishes on $N - \Z(T)$, and is a multiple of\/
$\zetaQ{T}$ on $\Z(T)$. Hence $\psi$ is the only character in $\Irr(\,N \mid \zetaQ{T}\,)$, and $G(\psi)$ is all of\/ $G$.
\end{proposition}
\begin{proof} Suppose that $\rho \in N - \Z(T)$. Then $\rho \Z(T)$ does not
lie in the radical $1$ of the symplectic $G/N$-group $\SSS(T) = N/\Z(T)$.
In view of \eqref{SCBar}, this gives us some $\sigma \in N$ such that
$\zetaQ{T}([\rho, \sigma]) = \bar c(\rho \Z(T), \sigma \Z(T)) \ne 1$.
Because $\zetaQ{T} = (\zetaQ{\psi})_{\Z(T)}$ by \eqref{ZetaQ}, it follows
from \eqref{PsiZetaPsi} that
\[  \psi(\rho) = \psi(\rho^{\sigma}) =\psi(\rho[\rho, \sigma]) =
\psi(\rho)\zetaQ{\psi}([\rho,\sigma]) = \psi(\rho)\zetaQ{T}([\rho,\sigma])
\in \CC. \]
Since $\zetaQ{T}([\rho,\sigma]) \ne 1$, this forces $\psi(\rho)$ to be $0$
for all $\rho \in N - \Z(T)$. 

Proposition \ref{UnderZetaQProp} tells us that $\psi$ is a multiple
$\psi(1)\zetaQ{T}$ of the $G$-invariant character $\zetaQ{T}$ on $\Z(T)$.
We have just seen that $\psi$ vanishes on $N - \Z(T)$. Hence the induced character $(\zetaQ{T})^N$ is a multiple $([N : \Z(T)]/\psi(1))\psi$ of $\psi$. This implies the rest of the proposition.
\end{proof}

In practice we apply the above considerations to various subtriples of $T$,
rather than to $T$ itself. Suppose that $\SSS(T')$ is abelian, for some
subtriple $T' = (G',N',\psi')$ of $T$. The natural embedding $e^T_{T'}$ in
\eqref{ETT} sends the factor group $G'/\N'$ isomorphically onto the
subgroup
\begin{subequations} \mylabel{ImETTEqs}
\begin{equation} \mylabel{ImETT} E^T_{T'} = e^T_{T'}(G'/N') = G'N/N
\end{equation}
of $G/N$. We use this isomorphism to translate the action \eqref{GNOnST} of
$G'/N'$ on $\SSS(T')$ to one of $E^T_{T'}$. Under this translated action,
the coset $\tau N = e^T_{T'}(\tau N') \in E^T_{T'}$ sends any
$\bar{\sigma} \in \SSS(T')$ to
\begin{equation} \mylabel{ImETTOnST} \bar{\sigma}^{\tau N} =
\bar{\sigma}^{\tau N'} = \bar{\sigma}^{\tau} \in \SSS(T'),
\end{equation}
\end{subequations}
for any $\tau \in G'$. Evidently the $G'/N'$-invariant function $\bar c =
\bar c_{T'} \colon \SSS(T') \times  \SSS(T') \to \CC^{\times}$  is also
$E^T_{T'}$-invariant. Thus $\SSS(T')$, with the bilinear form $\bar c$, is
a formed abelian $E^T_{T'}$-group.

Suppose that $A$ and $A'$ are formed abelian $H$-groups, with associated
bilinear forms $b_A$ and $b_{A'}$, for some finite group $H$. By an
\emph{isomorphism} $i$ of $A$ onto $A'$ as formed abelian $H$-groups we
mean an isomorphism $i$ of the group $A$ onto $A'$ which preserves the
actions of $H$ on those two groups, in the sense that
\begin{subequations} \mylabel{SSIsoEqs}
\begin{equation} \mylabel{IsoAct} i(\sigma^{\tau}) = i(\sigma)^{\tau} \in
A',
\end{equation}
for any $\sigma \in A$ and $\tau \in H$, and also preserves bilinear
forms, in the sense that
\begin{equation} \mylabel{IsoFun} b_{A'}(i(\rho), i(\sigma)) = b_A(\rho,
\sigma) \in \CC^{\times},
\end{equation}
\end{subequations}
for any $\rho, \sigma \in A$.  If such an $i$ exists, we say that the
formed abelian $H$-groups $A$ and $A'$ are \emph{isomorphic}. Clearly being
isomorphic is an equivalence relation among formed abelian $H$-groups.

We're going to show that equivalence among subtriples $T'$ of $T$ implies
isomorphism between their  associated formed abelian $E^T_{T'}$-groups. We
first treat the case of $N$-conjugate subtriples.
\begin{lemma} \mylabel{ConjIsoLemma} If\/ $\tau \in N$, then the section
$\SSS(T')$ is abelian, for some subtriple $T' = (G',N',\psi')$ of\/ $T$, if
and only if\/ $\SSS((T')^{\tau})$ is abelian. In that case both $G'/N'$ and
$(G')^{\tau}/(N')^{\tau}$ have the same image $E = E^T_{T'} =
E^T_{(T')^{\tau}}$ in $G/N$. Furthermore, conjugation by $\tau$ is an
isomorphism of\/ $\SSS(T')$ onto $\SSS((T')^{\tau})$ as formed abelian
$E$-groups.
\end{lemma}
\begin{proof}  Conjugation by $\tau \in N \le G(\psi)$ sends the normal
subgroup $N'$ in $T'$ to the normal subgroup $(N')^{\tau}$ in
$(T')^{\tau}$. It sends $\Z(T')$ onto $\Z((T')^{\tau}) = \Z(T')^{\tau}$ by
\eqref{ConjZT}. Hence it sends the factor group $\SSS(T') = N'/\Z(T')$
isomorphically onto $\SSS((T')^{\tau}) = (N')^{\tau}/\Z((T')^{\tau})$. In
particular, $\SSS(T')$ is abelian if and only if $\SSS((T')^{\tau})$ is
abelian.

Suppose that $\SSS(T')$ and $\SSS((T')^{\tau})$ are abelian. Since
$\tau$ lies in $N$, we have
\[  E^T_{(T')^{\tau}} = (G')^{\tau}N/N = (G'N)^{\tau}/N = G'N/N = E^T_{T'}.
\]
Thus $E = E^T_{T'} = E^T_{(T')^{\tau}}$ acts on both $\SSS(T')$ and
$\SSS((T')^{\tau})$ via \eqref{ImETTOnST}.  Any element of $E = G'N/N$
has the form $\sigma N$, for some $\sigma \in G'$. Because $\tau$ lies in
$N \normaleq G$, the conjugate $(\sigma N)^{\tau} = \sigma^{\tau} N$ is the
same element $\sigma N$ of $E$. It follows that
\[  (\bar{\rho}^{\sigma N})^{\tau} = \bar{\rho}^{\sigma\tau} =
\bar{\rho}^{\tau\sigma^{\tau}} = (\bar{\rho}^{\tau})^{\sigma^{\tau}N} =
(\bar{\rho}^{\tau})^{\sigma N}, \]
for any $\bar{\rho} \in \SSS(T')$.  So conjugation by $\tau$ is an
isomorphism of the group $\SSS(T')$ onto $\SSS((T')^{\tau})$, preserving
the actions of $E$ on these two abelian groups.

In view of \eqref{ConjZetaT}, \eqref{ImETTOnST}, and \eqref{SCBar} we have
\begin{multline*} \bar c_{(T')^{\tau}}((\rho \Z(T'))^{\tau}, (\sigma
\Z(T'))^{\tau}) = \bar c_{(T')^{\tau}}(\rho^{\tau} \Z((T')^{\tau}),
\sigma^{\tau}\Z((T')^{\tau})) = \zetaQ{(T')^{\tau}}([\rho^{\tau},
\sigma^{\tau}]) = \\
(\zetaQ{T'})^{\tau}([\rho,\sigma]^{\tau}) = \zetaQ{T'}([\rho,\sigma]) =
\bar c_{T'}(\rho \Z(T'), \sigma \Z(T')),  
\end{multline*}
for any $\rho, \sigma \in N'$. Therefore the $E$-isomorphism $\bar{\sigma}
\mapsto \bar{\sigma}^{\tau}$ of $\SSS(T')$ onto $\SSS((T')^{\tau})$
preserves  bilinear forms, and the lemma is proved.
\end{proof}

Next we consider covering subtriples.
\begin{lemma} \mylabel{CoverIsoLemma} If a subtriple $T' = (G',
N',\psi')$ covers a subtriple $T'' = (G'',N'',\psi'')$ of\/ $T$,
then $\SSS(T')$ is abelian if and only if\/ $\SSS(T'')$ is abelian. In that
case, $G'/N'$ and $G''/N''$ have the same image $E = E^T_{T'} = E^T_{T''}$
in $G/N$. Furthermore inclusion $G' \hookrightarrow G''$ induces an
isomorphism $i$ of\/ $\SSS(T')$ onto $\SSS(T'')$ as formed abelian
$E$-groups.
\end{lemma}
\begin{proof} We may apply all the results in \S \ref{Covers} with the
present $T'$ and $T''$ in place of the $\sq T$ and $T$ there. In
particular, \eqref{ITTGZs} tells us that inclusion $G' \hookrightarrow G''$
induces an isomorphism $i^{T''}_{T'}$ of $G'/\Z(T')$ onto $G''/\Z(T'') =
G'\Z(T'')/\Z(T'')$. This isomorphism sends $\SSS(T') = N'/\Z(T')$ onto
$\SSS(T'') = N''/\Z(T'') = N'\Z(T'')/\Z(T'')$, since $N'' \in \ZG(T'')$
corresponds to $N' = G' \cap N''$ in \eqref{CoverHCorr}. Hence $\SSS(T')$
is abelian if and only if $\SSS(T'')$  is abelian.

Now assume that $\SSS(T')$ and $\SSS(T'')$ are abelian. By \eqref{CoverG}
we have $G'' = G'\Z(T'')$. Since $\Z(T'') \le N'' = G'' \cap N$, it follows
that $E^T_{T''} = G''N/N = G'\Z(T'')N/N = G'N/N = E^T_{T'}$. 

We have already seen that the isomorphism $i^{T''}_{T'} \colon \bar{\sigma}
\mapsto \bar{\sigma}\Z(T'')$ of $G'/\Z(T')$ onto $G''/\Z(T'')$ restricts to
an isomorphism $i$ of $\SSS(T')$ onto $\SSS(T'')$. Any element
$\bar{\sigma}$ of $E = E^T_{T'} = G'N/N$ has the form $\sigma N$, for some
$\sigma \in G'$. Since $\sigma$ also lies in $G''$, and hence normalizes
$\Z(T'') \normaleq G''$, it follows from \eqref{ImETTOnST} that
\[  i(\bar{\rho}^{\bar{\sigma}}) = i(\bar{\rho}^{\sigma}) =
\bar{\rho}^{\sigma}\Z(T'') = (\bar{\rho}\Z(T''))^{\sigma} =
i(\bar{\sigma})^{\bar{\sigma}}, \]
for any $\bar{\rho} \in \SSS(T')$. Thus the isomorphism $i$ of $\SSS(T')$
onto $\SSS(T'')$ preserves the actions of $E$ on these two abelian groups.

We know from \eqref{CoverZeta} that $\zetaQ{T'}$ is the restriction of
$\zetaQ{T''}$ to $\Z(T') \le \Z(T'')$. It follows from this and
\eqref{SCBar} that
\[ \bar c_{T'}(\rho \Z(T'), \sigma\Z(T'))  = \zetaQ{T'}([\rho, \sigma]) =
\zetaQ{T''}([\rho,\sigma]) = \bar c_{T''}(i(\rho \Z(T')), i(\sigma \Z(T'))), \]
for any $\rho, \sigma \in N'$. Therefore $i$ preserves bilinear forms,
and the lemma is proved.
\end{proof}

Putting the above two lemmas together, we obtain
\begin{theorem} \mylabel{EquivIsoThm} If\/ $\SSS(T')$ is abelian, for
some subtriple $T'$ of\/ $T$, then so is $\SSS(T'')$, for any
subtriple $T''$ equivalent to $T'$. In that case both $E^T_{T'}$ and
$E^T_{T''}$ are the same subgroup $E$ of\/ $G/N$. Furthermore, the two
formed abelian $E$-groups $\SSS(T')$ and $\SSS(T'')$ are isomorphic.
\end{theorem}
\begin{proof} By Lemma \ref{ConjIsoLemma} this theorem is true when $T''$
is $N$-conjugate to $T'$. By Lemma \ref{CoverIsoLemma} it is true when $T'$
covers $T''$, or when $T''$ covers $T'$. An arbitrary triple $T''$
equivalent to $T'$ is obtained from $T'$ in a finite series of
equivalences, each of which is in one of these three cases. Since the
theorem holds at each step in this process, and its conclusions are clearly
transitive, this is enough to prove it in general.
\end{proof}

Now we can prove Theorem A of the introduction, in the form of
\begin{theorem} \mylabel{AThm} If\/ $\SSS(T')$ is nilpotent, for some
linear limit $T'$ of\/ $T  = (G,N,\psi) \in \mfrt$, then it is naturally a
$G(\psi)/N$-anisotropic symplectic $G(\psi)/N$-group. So is $\SSS(T'')$,
for any other linear limit $T''$ of\/ $T$. Furthermore, $\SSS(T')$ is
isomorphic to $\SSS(T'')$ as symplectic $G(\psi)/N$-groups.
\end{theorem}
\begin{proof} Since the linear limit $T'$ is linearly irreducible, and
$\SSS(T')$ is nilpotent, Propositions \ref{NilpSTProp} and \ref{AbIrrProp}
tell us that $\SSS(T')$ is abelian and $G'/N'$-anisotropic. So it is a
symplectic $G'/N'$-group. Now Proposition \ref{InvPsiProp} says that $G' =
G'(\psi')$. In view of Proposition \ref{MultLinRedStabProp}, this implies
that $G(\psi)/N = G'(\psi')N/N$ is the image $E^T_{T'} = G'N/N$ of
$G'/N'$ in $G/N$. So $\SSS(T')$ is a $G(\psi)/N$-anisotropic symplectic
$G(\psi)/N$-group, with the action \eqref{ImETTOnST} of $G(\psi)/N =
E^T_{T'}$ on $\SSS(T')$.

By Theorem \ref{LLThm} any linear limit $T''$ of $T$ is equivalent to $T'$.
The rest of the present theorem follows from this and Theorem
\ref{EquivIsoThm}.
\end{proof}

When $\SSS(T)$ is symplectic it is easy to classify the multilinear reductions of $T$.
\begin{proposition} \mylabel{SympMultRedProp} If\/ $\SSS(T)$ is symplectic, for some $T = (G,N,\psi) \in \mfrt$, then the map $f \colon T' \mapsto (\Z(T'), \zetaQ{T'})$ is a bijection of the set $\MLR(T)$ of all multilinear reductions of\/ $T$ onto the set $\mfrl = \mfrl(T)$ of all ordered pairs $(L,\lambda)$ satisfying
\[ L \normaleq G, \quad  \Z(T) \le L \le N \quad \text{and} \quad \lambda \in \Lin(\,L \mid \zetaQ{T}\,). \]
If\/ $(L,\lambda) \in \mfrl$, then the linear reduction $T(\lambda)$ of\/ $T$ is defined. The resulting map $g \colon (L,\lambda) \mapsto T(\lambda)$ from $\mfrl$ to $\MLR(T)$ is the inverse bijection to $f$. Hence any multilinear reduction of\/ $T$ is a linear reduction of\/ $T$.
\end{proposition}
\begin{proof}  Suppose that $T' = (G',N',\psi')$ is a multilinear reduction of $T$. We have $\Z(T) \le \Z(T') \le N' \le N$ by \eqref{MLRZZeta}.  Because $N/\Z(T) = \SSS(T)$ is abelian, this implies that $N$ normalizes $\Z(T')$. Propositions \ref{MultLinRedStabProp} and \ref{InvPsiProp} tell us that $G'(\psi')N = G(\psi) = G$. Since $G'(\psi')$ normalizes $\Z(T') \normaleq G'$, we conclude that $\Z(T')$ is a normal subgroup of $G$ such that $\Z(T) \le \Z(T') \le N$. The linear character $\zetaQ{T'}$ of $\Z(T')$ extends $\zetaQ{T} \in \Lin(\Z(T))$ by \eqref{MLRZZeta}. Therefore $(\Z(T'), \zetaQ{T'})$ belongs to $\mfrl$. So the function $f \colon \MLR(T) \to \mfrl$ in the proposition is defined.

If $(L,\lambda)$ belongs to $\mfrl$, then there is some $\phi \in \Irr(N)$ lying over $\lambda \in \Lin(L)$. Since $\phi$ lies over $\zetaQ{T} \le \lambda$, it must be the unique character $\psi \in \Irr(\,N \mid \zetaQ{T}\,)$ in Proposition \ref{InvPsiProp}. Hence $\lambda \le \psi$, and the linear reduction $T(\lambda)$ is defined. Thus the function $g \colon \mfrl \to \MLR(T)$ in the proposition is defined.

When $T' = (G',N',\psi')$ belongs to $\MLR(T)$, Proposition \ref{ZetaDetProp} tells us that $G' = G(\zetaQ{T'})$ and $N' = N(\zetaQ{T'})$. It also tells us that $\psi'$ is the unique character in $\Irr(\,N(\zetaQ{T'}) \mid \zetaQ{T'}\,)$ lying under $\psi$. Since we know that $\Z(T')$ is a normal subgroup of $N$, and that $\zetaQ{T'} \in \Lin(\Z(T'))$ lies under $\psi \in \Irr(N)$, this just says that $\psi'$ is the $\zetaQ{T'}$-Clifford correspondent $\psi_{\zetaQ{T'}}$ of $\psi$. Hence $T'$ is $T(\zetaQ{T'}) = (G(\zetaQ{T'}), N(\zetaQ{T'}), \psi_{\zetaQ{T'}})$. We conclude that the composite function $g \circ f$ is the identity map of $\MLR(T)$ onto itself.

Now fix $(L,\lambda) \in \mfrl$.  As in Step \ref{MLPerp} of the proof of Theorem \ref{LLThm}, the form $c \colon N \times N \to \CC^{\times}$ in \eqref{SC} satisfies
\[ c(\rho,\sigma) = \zetaQ{T}([\rho,\sigma]) = \lambda(\rho)\lambda^{\sigma^{-1}}(\rho) \]
for any $\rho \in L$ and $\sigma \in N$. It follows that the stabilizer $N(\lambda)$ is precisely the perpendicular subgroup $L^{\perp}$ to $L$ with respect to $c$. 

The commutator $[N(\lambda), \Z(T(\lambda))]$ is contained in $\Ker(\zetaQ{T(\lambda)})$, since $\zetaQ{T(\lambda)} \in \Lin(\Z(T(\lambda))$ is $N(\lambda)$-invariant. Because $\zetaQ{T(\lambda)}$ extends $\zetaQ{T}$, it follows that
\[  c(N(\lambda), \Z(T(\lambda))) = \zetaQ{T}([N(\lambda),\Z(T(\lambda))]) = \zetaQ{T(\lambda)}([N(\lambda),\Z(T(\lambda))]) = 1. \]
So $\Z(T(\lambda))$ is contained in $N(\lambda)^{\perp} = L^{\perp\perp}$. But $L^{\perp\perp}$ is the inverse image in $N$ of the double perpendicular subgroup $\bar L^{\perp\perp}$ to $\bar L = L/\Z(T)$ under the form $\bar c \colon (N/\Z(T)) \times (N/\Z(T)) \to \CC^{\times}$ induced by $c$. Since $\bar L^{\perp\perp} = \bar L$, by \eqref{BPerpPerp} for the symplectic $G/N$-group $\SSS(T) = N/\Z(T)$, we conclude that $\Z(T(\lambda)) \le L^{\perp\perp} = L$. The opposite inclusion holds by \eqref{LambdaUnderZeta}. Therefore $\Z(T(\lambda))$ is equal to $L$. Because $\lambda \in \Lin(L)$ is the restriction of $\zetaQ{T(\lambda)} \in \Z(T(\lambda))$  in \eqref{LambdaUnderZeta}, this implies that $\zetaQ{T(\lambda)} = \lambda$. Hence the composite map $f \circ g$ sends $(L,\lambda)$ to itself, and the proposition is proved.
\end{proof}
\begin{corollary} \mylabel{SympMultRedCor} If\/ $(L, \lambda)$ lies in $\mfrl$, then the normal subgroup $N(\lambda)$ in $T(\lambda)$ is the perpendicular subgroup $L^{\perp}$ to $L$ with respect to the form $c \colon N \times N \to \CC^{\times}$ in \eqref{SC}.
\end{corollary}
\begin{proof} This was shown in the course of the above proof.
\end{proof}

The linear limits of $T$ can also be described in the situation of the preceding proposition.
\begin{proposition} \mylabel{SympLinLimProp} If\/ $(K,\kappa)$ and $(L,\lambda)$ belong to the set $\mfrl$ in Proposition \ref{SympMultRedProp}, then $T(\kappa)$ is a linear reduction of\/ $T(\lambda)$ if and only if $L \le K$ and $\lambda = \kappa_L$.  Hence $T(\lambda)$ is a linear limit of\/ $T$ if and only if $\bar L = L/\Z(T)$ is maximal among all $G/N$-invariant isotropic subgroups of\/ $\SSS(T) = N/\Z(T)$.
\end{proposition}
\begin{proof} If $T(\kappa)$ is a linear reduction of $T(\lambda)$, then \eqref{MLRZZeta} tells us that $\Z(T(\lambda)) \le \Z(T(\kappa))$ and $\zetaQ{T(\lambda)} = (\zetaQ{T(\kappa)})_{\Z(T(\lambda))}$. In view of Proposition \ref{SympMultRedProp} this just says that $L \le K$ and $\lambda = \kappa_L$.

On the other hand, if $L \le K$ and $\lambda = \kappa_L$, then $T(\kappa)$ is the linear reduction $[T(\lambda)](\kappa)$ of $T(\lambda)$ by Proposition \ref{CompoundLinRedProp}. So the first statement in the present proposition is proved.

Since $L$ is a normal subgroup of $G$ contained in $N$, and $\lambda \in \Lin(L)$ lies under $\psi \in \Irr(N)$, Propositions \ref{IsotropicL} and \ref{AbIrrProp} tell us that $\bar L = L/\Z(T)$ is a $G/N$-invariant isotropic subgroup of $\SSS(T) = N/\Z(T)$.  
If $\bar L$ is maximal among such subgroups, then $\bar K = K/\Z(T)$, which is also $G/N$-invariant and isotropic, must equal $\bar L$. This forces $K$ to equal $L$, and $\kappa \in \Lin(K)$ to equal $\lambda = \kappa_L$.  So $T(\kappa) = T(\lambda)$. In view of Proposition \ref{SympMultRedProp}, this implies that $T(\lambda)$ is linearly irreducible. So $T(\lambda)$ is a linear limit of $T$ whenever $\bar L$ is maximal.

Suppose that there is some $G/N$-invariant isotropic subgroup $\bar M$ of $\SSS(T)$ such that $\bar L < \bar M$. Then $\bar M = M/\Z(T)$, where $M$ is a normal subgroup of $G$ satisfying $\Z(T) \le L < M \le N$.  Furthermore, $[M,M]$ is contained in $\Ker(T)$ by Proposition  \ref{AbIrrProp}. Because $\Ker(T) = \Ker(\zetaQ{T})$ and $\zetaQ{T} = \lambda_{\Z(T)}$, we conclude that $[M,M] \le \Ker(\lambda)$. So $\lambda \in \Lin(\,L \mid \zetaQ{T}\,)$ can be extended to some character $\mu \in \Lin(\,M \mid \zetaQ{T}\,)$. Then $(M,\mu)$ lies in $\mfrl$, and $T(\mu)$ is a linear reduction of $T(\lambda)$. This linear reduction is proper, since $M = \Z(T(\mu))$ strictly contains $L = \Z(T(\lambda))$. Therefore $T(\lambda)$ is linearly reducible when $\bar L$ is not maximal, and the remaining statement in the proposition is proved.
\end{proof}

Let $A$ be a formed abelian $H$-group, for some finite group $H$,
and $B$ be an isotropic $H$-invariant subgroup of $A$. Then
$B^{\perp}$ is an $H$-invariant subgroup of $A$ containing $B$, and the
factor group $B^{\perp}/B$ is naturally an $H$-group. The bilinear form $b
\colon A \times A \to \CC^{\times}$ for $A$ induces a bilinear form $\bar b
\colon (B^{\perp}/B) \times (B^{\perp}/B) \to \CC^{\times}$, sending the
cosets $\rho B, \sigma B \in B^{\perp}/B$ to
\begin{equation} \mylabel{BBar} \bar b(\rho B, \sigma B) = b(\rho, \sigma)
\in \CC^{\times},
\end{equation}
for any $\rho, \sigma \in B^{\perp}$. This induced form $\bar b$ is
strongly alternating and $H$-invariant, since $b$ has these properties. So
it turns $B^{\perp}/B$ into a formed abelian $H$-group. Whenever we treat
$B^{\perp}/B$ as a formed abelian $H$-group, it is this induced bilinear
form we have in mind.  

Any  subgroup of $B^{\perp}/B$ is the factor group
$C/B$, for some subgroup $C$ of $B^{\perp}$ containing $B$. The inclusions
\[ B \le B^{\perp\perp} \le C^{\perp} \le B^{\perp}, \]
in $A$, and the above definition of $\bar b$, imply that the perpendicular
subgroup $(C/B)^{\perp}$ to $C/B$ in $B^{\perp}/B$ is the image
\begin{subequations} \mylabel{BPerpBEqs}
\begin{equation} \mylabel{CBPerp} (C/B)^{\perp} = C^{\perp}/B
\end{equation}
of the perpendicular subgroup $C^{\perp}$ to $C$ in $A$. In particular, the
radical of $B^{\perp}/B$ is the image
\begin{equation}\mylabel{RadBPerpB} \rad(B^{\perp}/B) =
(B^{\perp}/B)^{\perp} = B^{\perp\perp}/B
\end{equation}
\end{subequations}
of $B^{\perp\perp}$.  When $A$ is symplectic this last image is $1$ by
\eqref{BPerpPerp}. Hence $B^{\perp}/B$ is naturally a symplectic $H$-group,
for any $H$-invariant isotropic subgroup $B$ of any symplectic $H$-group
$A$.

\medskip
We apply the above remarks when $A$ is the symplectic $G/N$-group $\SSS(T) = N/\Z(T)$ in Proposition \ref{SympMultRedProp}, and $B$ is its $G/N$-invariant isotropic subgroup $\bar L = L/\Z(T)$ in Proposition \ref{AbIrrProp}, for some $(L,\lambda) \in \mfrl$.  The normal subgroup $N(\lambda)$ in $T(\lambda)$ is $L^{\perp}$ by Corollary \ref{SympMultRedCor}. Hence its image in $\SSS(T)$ is
\begin{subequations} \mylabel{SympModNBarEqs}
\begin{equation} \mylabel{NLambdaBar}  N(\lambda)/\Z(T) = L^{\perp}/\Z(T) = (L/\Z(T))^{\perp} = \bar L^{\perp}.
\end{equation}
Since $L$ is $\Z(T(\lambda))$ by Proposition \ref{SympMultRedProp}, the factor group $\SSS(T(\lambda))$ is
\begin{equation} \mylabel{NBarTLambda} \SSS(T(\lambda)) = N(\lambda)/\Z(T(\lambda)) = L^{\perp}/L. 
\end{equation}
We conclude that there is a natural isomorphism $i = i_{T(\lambda)}$ of the group $\SSS(T(\lambda))$ onto $\bar L^{\perp}/\bar L$, sending
$\sigma \Z(T(\lambda)) = \sigma L \in \SSS(T(\lambda))$ to
\begin{equation} \mylabel{SympIDef} i(\sigma L) = (\sigma \Z(T))\bar L \in \bar L^{\perp}/\bar L = (L^{\perp}/\Z(T))/(L/\Z(T)),
\end{equation}
\end{subequations}
for any $\sigma \in N(\lambda) = L^{\perp}$.  Because $\bar L^{\perp}/\bar L$ is an abelian group, so is the isomorphic group $\SSS(T(\lambda))$. Hence $\SSS(T(\lambda))$ is naturally a formed abelian $G(\lambda)/N(\lambda)$-group.  By Propositions \ref{MultLinRedStabProp} and \ref{InvPsiProp} the embedding $G(\lambda) \hookrightarrow G$ induces an isomorphism of $G(\lambda)/N(\lambda)$ onto $G(\psi)/N = G/N$. We use this isomorphism to carry the conjugation action of $G(\lambda)/\N(\lambda)$ on $\SSS(T(\lambda))$ to an action of $G/N$ on $\SSS(T(\lambda))$. In this way $\SSS(T(\lambda))$ becomes a formed abelian $G/N$-group.
\begin{proposition} \mylabel{SympBarNTLambdaProp} If\/ $(L,\lambda)$ belongs to the set $\mfrl$ in Proposition \ref{SympMultRedProp}, then the map $i$ in \eqref{SympIDef} is an isomorphism of\/ $\SSS(T(\lambda))$ onto $\bar L^{\perp}/\bar L$ as formed abelian $G/N$-groups. Hence $\SSS(T(\lambda))$ is a symplectic $G(\lambda)/N(\lambda)$-group, as well as a symplectic $G/N$-group.
\end{proposition}
\begin{proof}  The natural alternating form $\bar c_{T(\lambda)}$ on $\SSS(T(\lambda)) \times \SSS(T(\lambda))$ is defined by \eqref{SCBar} with $T(\lambda)$ in place of $T$. Because $\zetaQ{T(\lambda)} = \lambda$ extends $\zetaQ{T}$, we have
\[
\bar c_{T(\lambda)}(\rho L, \sigma L) = \zetaQ{T(\lambda)}([\rho,\sigma]) = \zetaQ{T}([\rho,\sigma]) = \bar c_T(\rho \Z(T), \sigma \Z(T))
\]
for all $\rho, \sigma \in N(\lambda)$. So the isomorphism $i$ in \eqref{SympIDef} sends $\bar c_{T(\lambda)}$ to the form  $\Bar{\Bar c}$ on $(\bar L^{\perp}/\bar L) \times (\bar L^{\perp}/\bar L)$ induced by $\bar c = \bar c_T$ as in \eqref{BBar}. 

If $\tau \in G(\lambda)$ and $\sigma \in L^{\perp}$, then the element $(\sigma L)^{\tau N}$ of the $G/N$-group $\SSS(T(\lambda)) = L^{\perp}/L$ is $\sigma^{\tau} L$ by definition.  The isomorphism $i$ carries this to
\[  (\sigma^{\tau} \Z(T))\bar L = (\sigma \Z(T))^{\tau} \bar L = ((\sigma \Z(T))\bar L)^{\tau} = i(\sigma L)^{\tau}. \]
So $i$ preserves actions of $G/N$ as well as biliear forms, and thus is an isomorphism of formed abelian $G/N$-groups.

We know that $\bar L^{\perp}/\bar L$ is a symplectic $G/N$-group, since $\SSS(T)$ is one. Hence the isomorphic abelian $G/N$-group $\SSS(T(\lambda))$ is symplectic.  It follows that $\SSS(T(\lambda))$ is a symplectic $G(\lambda)/N(\lambda)$-group. So the proposition is proved.
\end{proof}

\end{document}